\documentclass[11pt,a4paper]{article}

\usepackage{amsmath,amssymb}
\usepackage{amsthm}
\usepackage[abbrev]{amsrefs}
\usepackage{latexsym}
\usepackage{bbm}
\usepackage{mathtools}
\usepackage{comment}
\usepackage{tikz}

\numberwithin{equation}{section}

\usepackage[
    a4paper,
    left=30mm,
    right=30mm,
    top=30mm,
    bottom=30mm
]{geometry}

\theoremstyle{definition}

\newtheorem{thm}{Theorem}[section]
\newtheorem{prop}[thm]{Proposition}
\newtheorem{lem}[thm]{Lemma}

\newtheorem{dfn}[thm]{Definition}

\theoremstyle{definition}

\newtheorem*{rem*}{Remark}

\newcommand{\R}{\mathbb{R}}

\newcommand{\Z}{\mathbb{Z}}

\newcommand{\Prob}{\mathbb{P}}
\newcommand{\E}{\mathbb{E}}

\newcommand{\ch}{\operatorname{ch}}
\newcommand{\dtv}{d_\mathrm{TV}}

\title{Cutoff for the random transposition-3-cycle shuffle}

\author{Soma Ishikawa}

\date{}

\begin{document}

\maketitle

%==================================================
% 概要
%==================================================

\begin{abstract}
In this paper, we introduce a new card-shuffling model called the random transposition-3-cycle shuffle, which is a modification of the random transposition shuffle. At each step, three positions $i$, $j$, and $k$ are chosen independently and uniformly at random from $\{1,2,\ldots,n\}$. If $i$, $j$, and $k$ are all distinct, the three cards at these positions are rearranged so that none of them remains in its original position. There are exactly two such rearrangements, corresponding to the two possible 3-cycles, each of which is chosen with equal probability. If exactly two of $i$, $j$, and $k$ coincide, the cards at the two distinct positions are swapped, corresponding to a transposition. If $i=j=k$, no change is made.

We prove that this shuffle exhibits cutoff at time $\frac13 n\log n$ with a window of order $n$ and explicitly determine its limit profile. Furthermore, we establish cutoff and determine the limit profile for a biased version of the shuffle, in which the relative probabilities of transpositions and 3-cycles are varied.
\end{abstract}

\begin{quote}
\small

\noindent\text{2020 Mathematics Subject Classification.}
Primary 60J10; Secondary 60B15.

\noindent\text{Keywords.}
Card shuffling, symmetric group, mixing time, cutoff phenomenon, limit profile. 
\end{quote}

\section{Introduction}

Let $i$ and $j$ be independent random variables uniformly distributed on $\{1,2,\ldots,n\}$. If $i\neq j$, we swap the cards at positions $i$ and $j$. If $i=j$, we do nothing. This shuffle is called the random transposition shuffle. More precisely, it can be formulated as a random walk on the symmetric group $S_n$ driven by the following probability distribution $T_2$:
\[
T_2(\sigma)=
\begin{cases}
\dfrac{1}{n},
    & \sigma=\mathrm{id},\\[1ex]
\dfrac{2}{n^2},
    & \sigma=(i\ j),\\[1ex]
0,
    & \text{otherwise}.
\end{cases}
\]

For two probability distributions $\mu$ and $\nu$ on a finite set $E$, the total variation distance is defined by
\[
\dtv(\mu,\nu)=\max_{A\subset E}|\mu(A)-\nu(A)|=\frac12\sum_{x\in E}|\mu(x)-\nu(x)|=\frac12d_1(\mu,\nu).
\]

Diaconis and Shahshahani (1981) \cite{paper1981} showed that the random transposition shuffle exhibits cutoff at time
\[
\frac12 n\log n
\]
with a window of size $O(n)$ as it converges to the uniform distribution $U_n$ on $S_n$. That is,
\[
\lim_{c\to-\infty}\liminf_{n\to\infty}\dtv\left(T_2^{*\left\lfloor\frac12 n\log n+cn\right\rfloor},U_n\right)=1,
\]
and
\[
\lim_{c\to\infty}\limsup_{n\to\infty}\dtv\left(T_2^{*\left\lfloor\frac12 n\log n+cn\right\rfloor},U_n\right)=0.
\]

Furthermore, Teyssier (2020) \cite{paper2020} showed that the limit profile of this shuffle is given by
\[
\lim_{n\to\infty}\dtv\left(T_2^{*\left\lfloor\frac12 n\log n+cn\right\rfloor},U_n\right)=\dtv\left(\operatorname{Poiss}(1+e^{-2c}),\operatorname{Poiss}(1)\right),\qquad c\in\R, 
\]
where $\operatorname{Poiss}(a)$ denotes the Poisson law of parameter $a$.
\vskip \baselineskip

In this paper, we consider a modification of the random transposition shuffle, which we call the random transposition-3-cycle shuffle. Let $i$, $j$, and $k$ be independent random variables uniformly distributed on $\{1,2,\ldots,n\}$. If $i$, $j$, and $k$ are all distinct, we rearrange the three cards at positions $i$, $j$, and $k$ so that none of them remains in its original position. There are two possible rearrangements, each of which is chosen with equal probability. If exactly two of $i$, $j$, and $k$ are equal, we swap the cards at the two distinct positions. For example, if $i=j\neq k$, we swap the cards at positions $i=j$ and $k$. If $i=j=k$, we do nothing.

More precisely, we define the shuffle as follows.

\begin{dfn}
Define a probability distribution $T_3$ on the symmetric group $S_n$ by
\[
T_3(\sigma)=
\begin{cases}
\dfrac{1}{n^2},
    & \sigma=\mathrm{id},\\[1ex]
\dfrac{6}{n^3},
    & \sigma=(i\ j),\\[1ex]
\dfrac{3}{n^3},
    & \sigma=(i\ j\ k),\\[1ex]
0,
    & \text{otherwise}.
\end{cases}
\]
The random walk on $S_n$ driven by $T_3$ is called the random transposition-3-cycle shuffle.
\end{dfn}

Our aim is to establish the cutoff phenomenon and determine the limit profile for the random transposition-3-cycle shuffle. Our first main result provides an upper bound on its mixing time. To prove this result, we follow an approach similar to that used by Diaconis and Shahshahani (1981)~\cite{paper1981} to establish the upper bound on the mixing time of the random transposition shuffle. A key ingredient in the proof is an estimate of the irreducible character ratios corresponding to 3-cycles. For this purpose, we use Theorem 4 of
Lassalle (2008)~\cite{paper2008}.

\begin{thm}\label{thm:upper-bound}
Let
\begin{equation}\label{k-1}
k=\left\lfloor\frac{1}{3}n\log n+cn \right\rfloor, \qquad c>0.
\end{equation}
Then there exists a constant $b>0$, independent of $c$, such that, for all sufficiently large $n$,
\[
\dtv(T_3^{*k}, U_n)\leq be^{-3c}.
\]
\end{thm}

Our second main result provides a matching lower bound on the mixing time of the random transposition-3-cycle shuffle. The proof is based on comparing $T_3^{*k}$ and $U_n$ on the set of permutations with at least $M$ fixed points, where $M$ is chosen to be approximately half the expected number of positions that have never been selected during the first $k$ steps of the shuffle.

\begin{thm}\label{thm:lower-bound}
Let
\[
k=\left\lfloor\frac{1}{3}n\log n-cn\right\rfloor, \qquad c>0.
\]
Then there exists a constant $b>0$, independent of $c$, such that, for all sufficiently large $n$,
\[
\dtv(T_3^{*k}, U_n)\geq1-be^{-3c}.
\]
\end{thm}

Theorems~\ref{thm:upper-bound} and~\ref{thm:lower-bound} imply that the random transposition-3-cycle shuffle exhibits cutoff at time
\[
\frac13 n\log n
\]
with a window of size $O(n)$. That is,
\[
\lim_{c\to-\infty}\liminf_{n\to\infty}\dtv\left(T_3^{*\left\lfloor\frac13 n\log n+cn\right\rfloor},U_n\right)=1,
\]
and
\[
\lim_{c\to\infty}\limsup_{n\to\infty}\dtv\left(T_3^{*\left\lfloor\frac13 n\log n+cn\right\rfloor},U_n\right)=0.
\]

Furthermore, we determine the limit profile of this shuffle. The proof follows an approach similar to that of
Teyssier (2020)~\cite{paper2020}. As in the proof of the upper bound, estimates of the irreducible character ratios corresponding to 3-cycles play a key role.

\begin{thm}\label{thm:limit-profile}
Let
\begin{equation}\label{k-2}
k=\left\lfloor\frac{1}{3}n\log n+cn\right\rfloor, \qquad c\in\R.
\end{equation}
Then
\[
\lim_{n\to\infty}\dtv(T_3^{*k}, U_n)=\dtv\left(\operatorname{Poiss}(1+e^{-3c}),\operatorname{Poiss}(1)\right).
\]
\end{thm}

Finally, we establish the existence of cutoff and determine the limit profile for a biased version of the random transposition-3-cycle shuffle, in which the probabilities of choosing transpositions and 3-cycles are varied. 

The remainder of this paper is organized as follows. In Section~2, we introduce some preliminary results. Theorems~\ref{thm:upper-bound}, \ref{thm:lower-bound}, and~\ref{thm:limit-profile} are proved in Sections~3, 4, and~5, respectively. Section~6 is devoted to the biased random transposition-3-cycle shuffle.
\vskip \baselineskip

\section{Preliminaries}

In this section, we introduce some tools for proving the cutoff phenomenon and the existence of a limit profile for the random transposition-3-cycle shuffle.

\subsection{Representation theory}

Let $G$ be a finite group and $\rho$ an irreducible representation of $G$. Denote the dimension and character of $\rho$ by $d_\rho$ and $\ch^\rho$, respectively. Let $\widehat{G}$ be the set of equivalence classes of irreducible representations of $G$, and $\widehat{G}^*$ the set of equivalence classes of nontrivial irreducible representations of $G$. For a probability distribution $P$ on $G$ that is constant on conjugacy classes, define
\[
\widehat{P}(\rho)=\sum_{g\in G}P(g)\rho(g),
\]
and set
\[
s_\rho=\frac{\operatorname{Tr}(\widehat{P}(\rho))}{d_\rho}.
\]

Diaconis and Shahshahani (1981) \cite{paper1981} used the following inequality to establish an upper bound for the random transposition shuffle.

\begin{lem}[{\cite[Lemma 14]{paper1981}}]\label{lem1981}
Let $P$ be a probability on $G$ constant on every conjugacy class, and let $U$ be the uniform distribution on $G$. Then, for every positive integer $t$, 
\[
d_1(P^{*t},U)^2\leq\sum_{\rho\in\widehat{G}^*}d_\rho^2|s_\rho|^{2t}.
\]
\end{lem}

Although this bound yields an upper bound on the distance to stationarity, it does not determine the limit profile. Teyssier (2020) \cite{paper2020} refined this approach as follows and determined the limit profile of the random transposition shuffle.

\begin{lem}[{\cite[Lemma 2.1]{paper2020}}]\label{lem2020}
Let $P$ be a probability on $G$ constant on every conjugacy class, and $U$ the uniform distribution on $G$. Then, for every positive integer $t$, 
\begin{equation}\label{equ2020}
d_1(P^{*t},U)=\frac{1}{|G|}\sum_{g\in G}\left|\sum_{\rho\in\widehat{G}^*}d_\rho s_\rho^t\overline{\ch^\rho(g)}\right|, 
\end{equation}
where $\overline{z}$ denotes the complex conjugate of $z$. 
\end{lem}

The following lemma is obtained by separating the right-hand side of (\ref{equ2020}) into a main term and an error term and estimating the error.

\begin{lem}[Approximation lemma~{\cite[Lemma 2.1]{paper2020}}]\label{d_1-2}
For any $S\subset\widehat G^*$,
\[
\left|d_1(P^{*t}, U)-\frac{1}{|G|}\sum_{g\in G}\left|\sum_{\rho\in S}d_\rho s_\rho^t\overline{\ch^\rho(g)}\right|\right|\leq\sum_{\rho\in\widehat{G}^*\setminus S}d_\rho|s_\rho|^t.
\]
\end{lem}

For the random transposition-3-cycle shuffle, we also use Lemma~\ref{lem1981} to obtain an upper bound and Lemma~\ref{d_1-2} to determine the limit profile.
\vskip \baselineskip

\subsection{Partitions and representations of the symmetric group}

For a positive integer $n$, a sequence of positive integers $\lambda=(\lambda_1,\lambda_2,\ldots,\lambda_m)$ is called a partition of $n$ if
\[
\lambda_1\geq\lambda_2\geq\cdots\geq\lambda_m>0, \qquad
\sum_{i=1}^m\lambda_i=n. 
\]
In this case, we write $\lambda\vdash n$.

We denote the Young diagram corresponding to a partition $\lambda$ by the same symbol $\lambda$. The box in row $i$ and column $j$ of the Young diagram $\lambda$ is denoted by $(i,j)$. We also denote the conjugate partition of $\lambda$ by $\lambda'$.

For a box $(i,j)$ in the Young diagram $\lambda$, define $c(i,j)=j-i$, which is called the content of the box $(i,j)$.

Let $\rho_\lambda$ be the irreducible representation corresponding to a partition $\lambda\vdash n$. Since the irreducible representations of $S_n$ are in one-to-one correspondence with the partitions of $n$, we shall identify an irreducible representation with its corresponding partition and sometimes denote the representation itself by $\lambda$. We write $\operatorname{ch}^\lambda$ and $d_\lambda$ for the character and dimension of $\rho_\lambda$, respectively. The following lemmas give properties of $d_\lambda$ and $d_{\lambda'}$. 
\begin{lem}[{\cite[Lemma 3]{paper1981}}]\label{d1}
It holds that
\[
\sum_{\rho\in\widehat{G}}d_\rho^2=|G|. 
\]
In particular, if $G=S_n$, then
\begin{equation}\label{d1-1}
\sum_{\lambda \vdash n}d_\lambda^2=n!.
\end{equation}
Consequently, for every $\lambda \vdash n$ we have
\begin{equation}\label{d1-2}
d_\lambda\leq\sqrt{n!}.
\end{equation}
\end{lem}

\begin{lem}[{\cite[Lemma 8]{paper1981}}]
Let $\lambda\vdash n$, and let $\lambda'$ be its conjugate partition. Then
\begin{equation}\label{lambda'}
d_\lambda=d_{\lambda'}
\end{equation}
\end{lem}

Using the hook-length formula (see Chapter 3 of M{\'e}liot (2017) \cite{book2017}), we obtain the following bounds on $d_\lambda$. 

\begin{lem}[{\cite[Proposition 3.1]{paper2020}}]\label{d2}
Let $\lambda=(\lambda_1,\lambda_2,\ldots,\lambda_m)\vdash n$ and $\lambda^*=(\lambda_2,\ldots,\lambda_m)$. Then
\begin{equation}\label{d2-1}
d_\lambda\leq\binom{n}{\lambda_1}d_{\lambda^*}. 
\end{equation}
In particular, by (\ref{d1-2}), 
\begin{equation}\label{d2-2}
d_\lambda\leq\binom{n}{\lambda_1}\sqrt{(n-\lambda_1)!}. 
\end{equation}
\end{lem}

Let $p(n)$ denote the number of partitions of $n$. We use the following bound on $p(n)$.

\begin{lem}[{\cite[equation (2.12)]{paper1981}}]\label{p(n)}
For every positive integer $n$,
\begin{equation}\label{p(n)-2}
p(n) \leq \exp\left(\pi\left(\frac{2}{3}n\right)^{\frac12}\right).
\end{equation}
\end{lem}
\vskip \baselineskip

\subsection{Character ratios and their estimates}

\begin{lem}[{\cite[Lemma 5]{paper1981}}]\label{lem1}
Let $\rho$ be an irreducible representation of a finite group $G$, and $P$ a probability on $G$ constant on every conjugacy class. For the $i$th conjugacy class, let $P_i$ denote the value of $P$, $n_i$ the number of elements in the class, and $\ch_i^\rho$ the value of the character $\ch^\rho$ on the class. Then, 
\[
\widehat{P}(\rho)=s_\rho I,\qquad s_\rho=\frac{1}{d_\rho}\sum_i P_i n_i \ch^\rho_i. 
\]
\end{lem}

Let $n$ be an integer and $\lambda$ a partition of $n$. Define the character ratios $r_2(\lambda)$ and $r_3(\lambda)$ by
\[
r_2(\lambda)=\frac{\ch^\lambda(1\ 2)}{d_\lambda},
\qquad
r_3(\lambda)=\frac{\ch^\lambda(1\ 2\ 3)}{d_\lambda}.
\]
Then, by Lemma~\ref{lem1} and the definition of $T_3$, we obtain
\begin{equation}\label{s_lambda}
s_\lambda=\frac{1}{n^2}+\frac{3(n-1)}{n^2}r_2(\lambda)+\frac{(n-1)(n-2)}{n^2}r_3(\lambda).
\end{equation}

By Lassalle (2008) \cite{paper2008}, $r_2(\lambda)$ and $r_3(\lambda)$ can be expressed as follows.

\begin{lem}[{\cite[Theorem 4]{paper2008}}]\label{r2, r3}
Let $\lambda\vdash n$. Then $r_2(\lambda)$ and $r_3(\lambda)$ are given respectively by
\begin{align*}
r_2(\lambda)&=\frac{2}{n(n-1)}\sum_{(i, j)\in \lambda}c(i,j), \\
r_3(\lambda)&=\frac{3}{n(n-1)(n-2)}\left( \sum_{(i, j)\in \lambda}{c(i,j)}^2-\binom{n}{2} \right). 
\end{align*}
\end{lem}

\begin{lem}
Let $\lambda\vdash n$, and let $\lambda'$ be its conjugate partition. Then
\begin{align}
r_2(\lambda)&=-r_2(\lambda'), \label{r2-lambda'}\\
r_3(\lambda)&=r_3(\lambda'). \label{r3-lambda'} 
\end{align}
\end{lem}
\begin{proof}
By the definition of conjugate partitions and Lemma~\ref{r2, r3}, these identities follow immediately. 
\end{proof}

We introduce the dominance order on partitions and establish the monotonicity of $r_2(\lambda)$ with respect to this order. This approach was also used by Diaconis and Shahshahani (1981)~\cite{paper1981}. 

\begin{dfn}
Let $\lambda=(\lambda_1,\ldots,\lambda_m)$ and $\mu=(\mu_1,\ldots,\mu_{m'})$ be partitions of $n$. If $m\leq m'$ and for every $1\leq k\leq m$
\[
\sum_{i=1}^{k}\lambda_i \geq \sum_{i=1}^{k}\mu_i, 
\]
then we write $\lambda \unrhd_2\mu$.
\end{dfn}

\begin{lem}[{\cite[Lemma 10]{paper1981}}]\label{r2}
Let $\lambda,\mu\vdash n$. If $\lambda\unrhd_2\mu$, then $r_2(\lambda) \geq r_2(\mu)$.
\end{lem}

To establish the monotonicity of $r_3(\lambda)$, we recall the Frobenius coordinates of partitions (see Macdonald (1995)~\cite[Chapter I, Section 1]{macdonald1995}) and introduce a relation on partitions based on these coordinates.

\begin{dfn}
Let $\lambda=(\lambda_1,\ldots,\lambda_m)\vdash n$ and let $\lambda'$ denote its conjugate partition. Let $d(\lambda)=\max\{i\geq1:\lambda_i\geq i\}$ denote the number of boxes on the main diagonal of the Young diagram $\lambda$. For each $1\leq i\leq d(\lambda)$, define $a_i(\lambda)=\lambda_i-i, \ b_i(\lambda)=\lambda_i'-i$. Then $\left(a_1(\lambda),\ldots,a_{d(\lambda)}(\lambda) \mid b_1(\lambda),\ldots,b_{d(\lambda)}(\lambda)\right)$ is called the Frobenius coordinates of $\lambda$.
\end{dfn}

\begin{dfn}
Let $\lambda,\mu\vdash n$. Let $\gamma_1(\lambda)\geq\gamma_2(\lambda)\geq\cdots$ be the non-increasing rearrangement of the $2d(\lambda)$ nonnegative integers $a_1(\lambda),\ldots,a_{d(\lambda)}(\lambda), \ b_1(\lambda),\ldots,b_{d(\lambda)}(\lambda)$ appearing in the Frobenius coordinates of $\lambda$. We also define $\gamma_k(\mu)$ in the same way. If $d(\lambda)\leq d(\mu)$ and for every $1\leq k\leq 2d(\lambda)$
\[
\sum_{i=1}^{k}\gamma_i(\lambda) \geq \sum_{i=1}^{k}\gamma_i(\mu), 
\]
then we write $\lambda\unrhd_3\mu$.
\end{dfn}

\begin{lem}\label{r3}
Let $\lambda,\mu\vdash n$. If $\lambda\unrhd_3\mu$, then $r_3(\lambda) \geq r_3(\mu)$.
\end{lem}

\begin{proof}
It suffices to show that $\sum_{(i,j)\in\lambda}(j-i)^2 \geq \sum_{(i,j)\in\mu}(j-i)^2$. Let the Frobenius coordinates of $\lambda$ be $(a_1,\ldots,a_{d(\lambda)} \mid b_1,\ldots,b_{d(\lambda)})$. For each $1\leq i\leq d(\lambda)$, the boxes to the right of the $i$th diagonal box have contents $1,2,\ldots,a_i$, while those below it have contents $-1,-2,\ldots,-b_i$. Hence,
\[
\sum_{(i,j)\in\lambda}(j-i)^2= \sum_{i=1}^{d(\lambda)}\left(\sum_{k=1}^{a_i}k^2+\sum_{k=1}^{b_i}k^2\right).
\]
Set
\[
F(r)=\sum_{k=1}^{r}k^2=\frac{r(r+1)(2r+1)}{6}, \qquad r\in\Z_{\geq0}.
\]
By the definition of $\gamma_i(\lambda)$, we have
\[
\sum_{(i,j)\in\lambda}(j-i)^2=\sum_{i=1}^{2d(\lambda)}F(\gamma_i(\lambda)).
\]
Similarly,
\[
\sum_{(i,j)\in\mu}(j-i)^2=\sum_{i=1}^{2d(\mu)}F(\gamma_i(\mu)).
\]

Next, set $x_i=\gamma_i(\lambda), y_i=\gamma_i(\mu)$. By appending zeros if necessary, we may assume that $(x_i)$ and $(y_i)$ have the same finite length $m$. Since $\lambda\unrhd_3\mu$, we have
\[
\sum_{i=1}^{k}x_i\geq \sum_{i=1}^{k}y_i, \qquad 1\leq k\leq m.
\]
For an integer $t\geq0$, define
\[
G_x(t)=\sum_{i=1}^{m}(x_i-t)_+, \qquad G_y(t)=\sum_{i=1}^{m}(y_i-t)_+,
\]
where $(z)_+=\max\{z,0\}$. Since $(x_i)$ and $(y_i)$ are arranged in nonincreasing order, we have
\[
G_x(t)=\max_{0\leq k\leq m}\left\{\sum_{i=1}^{k}x_i-kt\right\}, \qquad
G_y(t)=\max_{0\leq k\leq m}\left\{\sum_{i=1}^{k}y_i-kt\right\}.
\]
Therefore, $\lambda\unrhd_3\mu$ implies that
\[
G_x(t)\geq G_y(t), \qquad \text{for all integers $t\geq0$}.
\]
For every nonnegative integer $r$,
\[
F(r)=r+\sum_{t=1}^{\infty}(2t+1)(r-t)_+.
\]
Consequently,
\begin{align*}
\sum_{i=1}^mF(x_i)&=\sum_{i=1}^m x_i+\sum_{t=1}^{\infty}(2t+1)G_x(t),\\
\sum_{i=1}^mF(y_i)&=\sum_{i=1}^m y_i+\sum_{t=1}^{\infty}(2t+1)G_y(t).
\end{align*}
Since $\lambda\unrhd_3\mu$, we have
\[
\sum_{i=1}^m x_i \geq \sum_{i=1}^m y_i.
\]
Moreover, since $G_x(t)\geq G_y(t)$ for every integer $t\geq0$ and $2t+1>0$, we obtain
\[
\sum_{i=1}^mF(x_i) \geq \sum_{i=1}^mF(y_i).
\]
Therefore,
\[
\sum_{(i,j)\in\lambda}(j-i)^2 \geq \sum_{(i,j)\in\mu}(j-i)^2.
\]
This completes the proof. 
\end{proof}

In Lemmas~\ref{r2-1}, \ref{r2-2}, \ref{r2-3}, \ref{r2-4}, and \ref{r2-5}, we estimate $r_2(\lambda)$ by considering different ranges of $\lambda_1$. The statements and proofs of these lemmas are based on Lemmas 11, 12, and 13 of Diaconis and Shahshahani (1981) \cite{paper1981}, with some additions and modifications made for the present paper.

\begin{lem}\label{r2-1}
Let $\lambda=(\lambda_1,\ldots,\lambda_m)\vdash n$ satisfy $\lambda_1\leq n/3$. Then
\[
r_2(\lambda) \leq \frac{1}{n(n-1)}\left(\frac{n^2}{3}-3n+\frac{20}{3}\right).
\]
\end{lem}

\begin{proof}
Let $b=\left\lceil n/3 \right\rceil$, and write
\[
b=\frac{n}{3}+\epsilon, \qquad 0\leq\epsilon\leq\frac23.
\]
Since $(b,b,n-2b)\unrhd_2\lambda$, Lemma~\ref{r2} implies that
\begin{align*}
r_2(\lambda)& \leq r_2((b,b,n-2b))\\
&=\frac{2}{n(n-1)}\left(\sum_{i=0}^{b-1}i+\sum_{i=-1}^{b-2}i+\sum_{i=-2}^{n-2b-3}i\right)\\
&=\frac{1}{n(n-1)}\left(b(b-1)+(b-2)(b-1)-2+(n-2b-3)(n-2b-2)-6\right)\\
&=\frac{1}{n(n-1)}\left(\frac{n^2}{3}-3n+6\epsilon^2+6\epsilon\right), 
\end{align*}
from which the conclusion follows readily.
\end{proof}

Arguing as in the proof of Lemma~\ref{r2-1}, we obtain the following five lemmas.

\begin{lem}\label{r2-2}
Let $\lambda=(\lambda_1,\ldots,\lambda_m)\vdash n$ satisfy $\lambda_1\leq n/2$. Then
\[
r_2(\lambda) \leq \frac{1}{n(n-1)}\left(\frac{n^2}{2}-2n+\frac32\right).
\]
\end{lem}

\begin{lem}\label{r2-3}
Let $\lambda=(\lambda_1,\ldots,\lambda_m)\vdash n$ satisfy $\lambda_1\leq 2n/3$. Then
\[
r_2(\lambda) \leq \frac{1}{n(n-1)}\left(\frac{5}{9}n^2-\frac{11}{9}n+\frac{20}{9}\right).
\]
\end{lem}

\begin{lem}\label{r2-4}
Let $\lambda=(\lambda_1,\ldots,\lambda_m)\vdash n$ satisfy $\lambda_1=n-j\geq 2n/3$. Then
\[
r_2(\lambda) \leq 1-\frac{2j(n-j+1)}{n(n-1)}.
\]
\end{lem}

\begin{lem}\label{r2-6}
Let $\lambda=(\lambda_1,\ldots,\lambda_m)\vdash n$ satisfy $\lambda_1>n/2$. Then
\[
r_2(\lambda) \geq 0.
\]
\end{lem}

\begin{lem}\label{r2-5}
Let $\lambda=(\lambda_1,\ldots,\lambda_m)\vdash n$ satisfy $\lambda_1\geq 2n/3$. Then
\[
r_2(\lambda) \geq \frac{n-3}{3(n-1)}.
\]
\end{lem}

In Lemmas~\ref{r3-1}, \ref{r3-2}, \ref{r3-3}, \ref{r3-4}, and \ref{r3-5}, we estimate $r_3(\lambda)$ by considering several cases
depending on the values of $\lambda_1$ and $m$.

\begin{lem}\label{r3-1}
Let $\lambda=(\lambda_1,\ldots,\lambda_m)\vdash n$ satisfy $\lambda_1\leq n/3$ and $m\leq n/3$. Then
\[
r_3(\lambda) \leq \frac{1}{n(n-1)(n-2)}\left(\frac{n^3}{9}-\frac{7}{3}n^2+\frac{40}{3}n+\frac{224}{3}\right).
\]
\end{lem}

\begin{proof}
Let $b=\left\lfloor n/3 \right\rfloor$, and write
\[
n=3b+r, \qquad r=0,1,2.
\]
Writing $\lambda$ in Frobenius coordinates as $\lambda=(a_1,\ldots,a_{d(\lambda)}\mid b_1,\ldots,b_{d(\lambda)})$, we see that the assumptions imply
\[
a_1,b_1 \leq b-1, \qquad a_2,b_2 \leq b-2 .
\]
Since $(b-1,b-2\mid b-1,r+2)\unrhd_3\lambda$, Lemma~\ref{r3} implies that
\begin{align*}
r_3(\lambda) &\leq r_3((b-1,b-2 \mid b-1,r+2))\\
&=\frac{3}{n(n-1)(n-2)}\left(2\sum_{i=1}^{b-1}i^2+\sum_{i=1}^{b-2}i^2+\sum_{i=1}^{r+2}i^2-\binom{n}{2}\right)\\
&=\frac{1}{n(n-1)(n-2)}(3b^3-21b^2-9br+12b+r^3+6r^2+20r+12)\\
&\leq \frac{1}{n(n-1)(n-2)}\left(\frac{n^3}{9}-\frac{7}{3}n^2+\frac{40}{3}n+\frac{224}{3}\right), 
\end{align*}
where the last inequality follows by substituting $b=(n-r)/3$ and considering the three cases $r=0,1,2$.
\end{proof}

Arguing as in the proof of Lemma~\ref{r3-1}, we obtain the following three lemmas.

\begin{lem}\label{r3-2}
Let $\lambda=(\lambda_1,\ldots,\lambda_m)\vdash n$ satisfy $\lambda_1\leq n/2$ and $m\leq n/2$. Then
\[
r_3(\lambda) \leq \frac{n-8}{4(n-2)}.
\]
\end{lem}

\begin{lem}\label{r3-3}
Let $\lambda=(\lambda_1,\ldots,\lambda_m)\vdash n$ satisfy $\lambda_1\leq 2n/3$ and $m\leq 2n/3$. Then
\[
r_3(\lambda) \leq \frac{n-2}{3(n-1)}.
\]
\end{lem}

\begin{lem}\label{r3-4}
Let $\lambda=(\lambda_1,\ldots,\lambda_m)\vdash n$ satisfy $\lambda_1=n-j\geq 2n/3$. Then
\[
r_3(\lambda) \leq 1-\frac{3j(n-j-1)}{(n-1)(n-2)}.
\]
\end{lem}

We also use the following lemma in the proof of Theorem~\ref{thm:upper-bound}. 

\begin{lem}\label{r3-5}
For every $\lambda=(\lambda_1,\ldots,\lambda_m)\vdash n$,
\[
r_3(\lambda) \geq -\frac{3}{2(n-2)}.
\]
In particular, if $\lambda_1\geq n/3$ or $m\geq n/3$, then for all sufficiently large $n$, 
\[
r_3(\lambda)\geq0. 
\]
\end{lem}
\begin{proof}
Since the sum of squared contents is nonnegative, we have
\begin{align*}
r_3(\lambda)&\geq -\frac{3}{n(n-1)(n-2)}\cdot\frac{n(n-1)}{2}\\
&=-\frac{3}{2(n-2)}.
\end{align*}
Furthermore, if $\lambda_1\geq n/3$ or $m\geq n/3$, then
\begin{align*}
r_3(\lambda)&\geq\frac{3}{n(n-1)(n-2)}\left(\sum_{i=0}^{\left\lceil\frac{n}{3}\right\rceil-1}i^2-\binom{n}{2}\right)\\
&\geq\frac{3}{n(n-1)(n-2)}\cdot\frac{n(n^2-45n+45)}{81}.
\end{align*}
Therefore, for all sufficiently large $n$, $r_3(\lambda)\geq0$.
\end{proof}
\vskip \baselineskip

\section{Upper bound (Proof of Theorem \ref{thm:upper-bound})}

In this section, we prove Theorem~\ref{thm:upper-bound}. Before proving the theorem, we establish the following lemma. For a positive integer $n$, let
\[
S(n)=\sum_{j=1}^{\left\lfloor\frac{n}{3}\right\rfloor}\frac{p(j)}{j!}\exp\left(\frac{2\log n}{n}j^2\right). 
\]

\begin{lem}
It holds that
\begin{equation}\label{S(n)}
\sup_{n\geq1}S(n)<\infty. 
\end{equation}
\end{lem}

\begin{proof}
Write
\begin{align*}
S(n)&=\sum_{j=1}^{\left\lfloor\frac{n}{\log n}\right\rfloor}\frac{p(j)}{j!}\exp\left(\frac{2\log n}{n}j^2\right)+\sum_{j=\left\lfloor\frac{n}{\log n}\right\rfloor+1}^{\left\lfloor\frac{n}{3}\right\rfloor}\frac{p(j)}{j!}\exp\left(\frac{2\log n}{n}j^2\right)\\
&=S_1(n)+S_2(n), 
\end{align*}
where $S_1(n)$ and $S_2(n)$ denote the first and second sums, respectively.

For $S_1(n)$, since $j\leq n/\log n$, we have
\[
\frac{2\log n}{n}j^2\leq2j.
\]
By (\ref{p(n)-2}) and $j!\geq(j/e)^j$, we obtain
\[
\frac{p(j)}{j!}\exp\left(\frac{2\log n}{n}j^2\right)
\leq\exp\left(\pi\left(\frac23j\right)^{\frac12}+3j-j\log j\right).
\]
Therefore,
\[
\sup_{n\geq1}S_1(n)<\infty. 
\]
For $S_2(n)$, since $j\leq n/3$, we have
\[
\frac{2\log n}{n}j^2 \leq \frac23j\log n.
\]
Also, since $j>n/\log n$ and $\log j>\log n-\log\log n$,
\begin{align*}
\frac{p(j)}{j!}\exp\left(\frac{2\log n}{n}j^2\right)
&\leq\exp\left(\pi\left(\frac23j\right)^{\frac12}-j\log j+j+\frac23j\log n\right)\\
&\leq\exp\left(-\frac13j\log n+j\log\log n+j+\pi\left(\frac23j\right)^{\frac12}\right). 
\end{align*}
For all sufficiently large $n$, $\log\log n+1+\pi\sqrt{\dfrac{2}{3j}}\leq\dfrac16\log n$. Therefore,
\[
\frac{p(j)}{j!}\exp\left(\frac{2\log n}{n}j^2\right)\leq\exp\left(-\frac16j\log n\right)\leq e^{-\frac{n}{6}}.
\]
Consequently,
\[
S_2(n)\leq\frac n3e^{-\frac{n}{6}}\xrightarrow[n\to\infty]{}0.
\]
Combining these estimates for $S_1(n)$ and $S_2(n)$, we obtain (\ref{S(n)}). 
\end{proof}

\begin{proof}[\normalfont\bfseries Proof of Theorem \ref{thm:upper-bound}]
Applying Lemma~\ref{lem1981} to $T_3$, we obtain
\[
d_1(T_3^{*k},U_n)^2\leq\sum_{\lambda\in\widehat{S_n}^*}d_\lambda^2|s_\lambda|^{2k}.
\]
We estimate the sum on the right-hand side by considering the following zones:
\begin{align*}
&\text{Zone 1. }\qquad  1 \leq \lambda_1 \leq \frac{n}{3} \  \text{and} \  1 \leq m \leq \frac{n}{3}; \\
&\text{Zone 2. }\qquad  \frac{n}{3}<\max\{\lambda_1,m\}<\frac{2}{3}n; \\
&\text{Zone 3. }\qquad  (1)\ \frac{2}{3}n \leq \lambda_1 \leq n-1 \qquad (2)\ \frac{2}{3}n \leq m \leq n-1 \qquad (3)\ m=n. \\
\end{align*}

\begin{center}
\begin{tikzpicture}[scale=4, xscale=1.5]
\draw[->] (0,0) -- (1.1,0)
  node[right] {$\lambda_1$};

\draw[->] (0,0) -- (0,1.1)
  node[above] {$m$};

\draw (0.05,0.05) -- (0.05,1);
\draw (1/3,0.05) -- (1/3,1/3);
\draw (2/3,0.05) -- (2/3,1/3+0.05);
\draw (0.05,0.05) -- (1,0.05);
\draw (0.05,1/3) -- (1/3,1/3);
\draw (0.05,2/3) -- (1/3+0.05,2/3);
\draw (0.05,1) -- (1,0.05);

\node[below] at (0.05,0) {$1$};
\node[below] at (1/3,0) {$n/3$};
\node[below] at (2/3,0) {$2n/3$};
\node[below] at (1,0) {$n$};
\node[left] at (0,0.05) {$1$};
\node[left] at (0,1/3) {$n/3$};
\node[left] at (0,2/3) {$2n/3$};
\node[left] at (0,1.) {$n$};

\draw[->] (0.2,1.05) -- (0.05,1);
\node[right] at (0.2,1.05) {3-(3)};

\node at (1/6+0.025,1/6+0.025) {1};
\node at (5/12,5/12) {2};
\node at (7/9,1/12+0.025) {3-(1)};
\node at (1/9+0.025,3/4) {3-(2)};
\end{tikzpicture}
\end{center}

\noindent\textbf{Zone 1}

By Lemmas~\ref{r2-1}, \ref{r3-1}, \ref{r3-5}, (\ref{s_lambda}) and (\ref{r2-lambda'}), we have
\begin{align*}
s_\lambda
&\leq \frac{1}{n^2}+\frac{3}{n^3}\left(\frac{n^2}{3}-3n+\frac{20}{3}\right)+\frac{1}{n^3}\left(\frac{n^3}{9}-\frac{7}{3}n^2+\frac{40}{3}n+\frac{224}{3}\right)\\
&=\frac{1}{9}-\frac{4}{3n}+\frac{16}{3n^2}+\frac{284}{3n^3}\\
s_\lambda& \geq \frac{1}{n^2}-\frac{3}{n^3}\left(\frac{n^2}{3}-3n+\frac{20}{3}\right)-\frac{3(n-1)}{2n^2}\\
&=-\frac{5}{2n}+\frac{23}{2n^2}-\frac{20}{n^3}.
\end{align*}
Thus, $|s_\lambda|\leq 1/9$ for all sufficiently large $n$. By the definition (\ref{k-1}) of $k$, 
\[
k\geq\frac{1}{3}n\log n+cn-1. 
\]
Therefore, using (\ref{d1-1}) and $n!\leq n^n$, we obtain
\begin{align*}
\sum_{\lambda:\text{Zone 1}}d_\lambda^2s_\lambda^{2k}
&\leq \left(\frac19\right)^{2k}n^n\\
&\leq\exp \left(\left(1-\frac43\log 3\right)n\log n-4cn\log3+4\log3\right)\\
&\xrightarrow[n\to\infty]{}0.
\end{align*}

\noindent\textbf{Zone 2}

Arguing as in Zone 1 and applying Lemmas~\ref{r2-3}, \ref{r3-3}, \ref{r3-5}, (\ref{s_lambda}) and (\ref{r2-lambda'}), we obtain
\[
-\frac{19}{6n}+\frac{37}{6n^2}-\frac{20}{3n^3} \leq s_\lambda \leq \frac{1}{3}+\frac{1}{3n}-\frac{4}{3n^2}+\frac{20}{3n^3}. 
\]
Hence, $|s_\lambda|\leq\dfrac13+o(1)$. (\ref{d2-2}) implies $d_\lambda\leq 2^n n^{n/3}$, so
\begin{align*}
\sum_{\lambda:\text{Zone 2}}d_\lambda^2s_\lambda^{2k}
&\leq p(n)4^n n^{\frac23n}\left(\frac13+o(1)\right)^{2k}\\
&\leq\exp \left(\frac23\left(1+\log\left(\frac13+o(1)\right)\right)n\log n+O(n)\right)\\
&\xrightarrow[n\to\infty]{}0.
\end{align*}

\noindent\textbf{Zone 3-(1)}

Write $\lambda_1=n-j$. By Lemmas~\ref{r2-4}, \ref{r2-6}, \ref{r3-4}, \ref{r3-5} and (\ref{s_lambda}), we have
\begin{align*}
0\leq s_\lambda 
&\leq1-\frac{6j(n-j+1)}{n^3}-\frac{3j(n-j-1)}{n^2}\\
&\leq 1-3\left(\frac{j}{n}-\frac{j^2}{n^2}+\frac{1}{n^2}\right). 
\end{align*}
Hence, for all sufficiently large $n$,
\[
|s_\lambda|\leq 1-3\left(\frac{j}{n}-\frac{j^2}{n^2}+\frac{1}{n^2}\right)\leq \exp\left(-3\left(\frac{j}{n}-\frac{j^2}{n^2}+\frac{1}{n^2}\right)\right).
\]
By (\ref{d2-2}), 
\begin{equation}\label{d3}
d_\lambda \leq \dfrac{n^j}{j!}\sqrt{j!}=\dfrac{n^j}{\sqrt{j!}}. 
\end{equation}
Therefore, 
\begin{align*}
&\sum_{\lambda:\text{Zone 3-(1)}}d_\lambda^2s_\lambda^{2k}\\
&=\sum_{j=1}^{\left\lfloor\frac{n}{3}\right\rfloor}\sum_{\lambda_1=n-j}d_\lambda^2s_\lambda^{2k}\\
&\leq \sum_{j=1}^{\left\lfloor\frac{n}{3}\right\rfloor}p(j)\frac{n^{2j}}{j!}\exp\left(-6k\left(\frac{j}{n}-\frac{j^2}{n^2}+\frac{1}{n^2}\right)\right)\\
&\leq\sum_{j=1}^{\left\lfloor\frac{n}{3}\right\rfloor}\frac{p(j)}{j!}\exp\left(\frac{2\log n}{n}(j^2-1)-6cn\left(\frac{j}{n}-\frac{j^2}{n^2}+\frac{1}{n^2}\right)+6\left(\frac{j}{n}-\frac{j^2}{n^2}+\frac{1}{n^2}\right)\right).
\end{align*}
Since $1\leq j\leq n/3$, we have
\[
\frac1n\leq\frac{j}{n}-\frac{j^2}{n^2}+\frac{1}{n^2}\leq\frac13
\]
and hence
\[
\sum_{\lambda:\text{Zone 3-(1)}}d_\lambda^2s_\lambda^{2k}
\leq e^{-6c+2}\sum_{j=1}^{\left\lfloor\frac{n}{3}\right\rfloor}\frac{p(j)}{j!}\exp\left(\frac{2\log n}{n}(j^2-1)\right).
\]
By (\ref{S(n)}),
\[
\sup_{n\geq1}\sum_{j=1}^{\left\lfloor\frac{n}{3}\right\rfloor}\frac{p(j)}{j!}\exp\left(\frac{2\log n}{n}(j^2-1)\right)\leq\sup_{n\geq1}S(n)<\infty. 
\]
Thus, there exists a constant $B>0$, independent of $c$, such that, for all sufficiently large $n$, 
\[
\sum_{\lambda:\text{Zone 3-(1)}}d_\lambda^2s_\lambda^{2k} \leq Be^{-6c}.
\]

\noindent\textbf{Zone 3-(2)}

Let $\lambda'$ be the conjugate partition of $\lambda$. By (\ref{s_lambda}), (\ref{r2-lambda'}) and (\ref{r3-lambda'}), we obtain
\begin{align*}
s_\lambda&=\frac{1}{n^2}-\frac{3(n-1)}{n^2}r_2(\lambda')+\frac{(n-1)(n-2)}{n^2}r_3(\lambda')\\
&=s_{\lambda'}-\frac{6(n-1)}{n^2}r_2(\lambda').
\end{align*}
Since $2n/3\leq\lambda'_1\leq n-1$, Lemma~\ref{r2-5} gives
\[
s_\lambda \leq s_{\lambda'}-\frac{2(n-3)}{n^2}.
\]
If $s_\lambda\geq0$, then
\[
|s_\lambda|\leq\left(1-\frac{2(n-3)}{n^2}\right)|s_{\lambda'}|.
\]
If $s_\lambda<0$, then $\lambda'$ belongs to Zone~3-(1), so $s_{\lambda'}\geq0$. Hence, 
\begin{align*}
|s_\lambda|=-s_\lambda=\frac{6(n-1)}{n^2}r_2(\lambda')-s_{\lambda'}
\leq\frac{6(n-1)}{n^2}\leq\frac6n. 
\end{align*}
Thus, 
\[
s_\lambda^{2k}\leq\left(1-\frac{2(n-3)}{n^2}\right)^{2k}s_{\lambda'}^{2k}+\left(\frac6n\right)^{2k}.
\]
Therefore, by (\ref{lambda'}), 
\begin{align*}
\sum_{\lambda:\text{Zone 3-(2)}}d_{\lambda}^2s_\lambda^{2k}
&\leq \sum_{\lambda:\text{Zone 3-(2)}}d_{\lambda}^2\left(\left(1-\frac{2(n-3)}{n^2}\right)^{2k}s_{\lambda'}^{2k}+\left(\frac6n\right)^{2k}\right)\\
&\leq \left(\sum_{\lambda:\text{Zone 3-(1)}}d_\lambda^2s_\lambda^{2k}\right)\left(1-\frac{2(n-3)}{n^2}\right)^{2k}+\sum_{\lambda:\text{Zone 3-(1)}}d_{\lambda}^2\left(\frac6n\right)^{2k}. \\
\end{align*}
For the first term,
\begin{align*}
&\left(\sum_{\lambda:\text{Zone 3-(1)}}d_\lambda^2s_\lambda^{2k}\right)\left(1-\frac{2(n-3)}{n^2}\right)^{2k}\\
&\leq\left(\sum_{\lambda:\text{Zone 3-(1)}}d_\lambda^2s_\lambda^{2k}\right)\exp\left(-2\left(\frac{1}{3}n\log n+cn-1\right)\left(\frac{2}{n}-\frac{6}{n^2}\right)\right)\\
&=\left(\sum_{\lambda:\text{Zone 3-(1)}}d_\lambda^2s_\lambda^{2k}\right)\exp\left(-\frac43\log n-4c+O\left(\frac{\log n}{n}\right)\right)\\
&\xrightarrow[n\to\infty]{}0.
\end{align*}
For the second term, 
\begin{align*}
\sum_{\lambda:\text{Zone 3-(1)}}d_{\lambda}^2\left(\frac6n\right)^{2k}
&\leq n^n\left(\frac6n\right)^{2k}\\
&=\exp\left(-\frac23n(\log n)^2+O(n\log n)\right)\\
&\xrightarrow[n\to\infty]{}0.
\end{align*}
Hence,
\[
\sum_{\lambda:\text{Zone 3-(2)}}d_{\lambda}^2s_\lambda^{2k}\xrightarrow[n\to\infty]{}0.
\]

\noindent\textbf{Zone 3-(3)}

In this case, $\lambda=(1,\ldots,1)$, which we denote by $(1^n)$. Since $r_2((1^n))=-1$ and $r_3((1^n))=1$, (\ref{s_lambda}) gives
\begin{align*}
s_{(1^n)}&=\frac{1}{n^2}-\frac{3(n-1)}{n^2}+\frac{(n-1)(n-2)}{n^2}\\
&=1-\frac6n+\frac{6}{n^2}.
\end{align*}
Since $d_{(1^n)}=1$,
\begin{align*}
d_{(1^n)}^2s_{(1^n)}^{2k}
&=\left(1-\frac6n+\frac{6}{n^2}\right)^{2k}\\
&\leq \exp \left(2\left(-\frac6n+\frac{6}{n^2}\right)\left(\frac13n\log n+cn-1\right)\right)\\
&=\exp \left(-4\log n-12c+O\left(\frac{\log n}{n}\right)\right)\\
&\xrightarrow[n\to\infty]{}0.
\end{align*}
\vskip \baselineskip

Combining these estimates, we conclude that there exists a constant $B'>0$, independent of $c$, such that, for all sufficiently large $n$, 
\[
d_1(T_3^{*k}, U_n)^2 \leq \sum_{\lambda:\text{Zone 1-3}}d_{\lambda}^2s_\lambda^{2k} \leq B'e^{-6c}.
\]
The desired result follows from this.
\end{proof}
\vskip \baselineskip

\section{Lower bound (Proof of Theorem \ref{thm:lower-bound})}

\begin{proof}[\normalfont\bfseries Proof of Theorem \ref{thm:lower-bound}]
Let $Y$ denote the number of positions that have never been selected among the $3k$ card selections made up to time $k$. For each $i\in\{1,\ldots,n\}$, let $I_i$ be the indicator that position $i$ has never been selected by time $k$. Then $Y=\sum_{i=1}^n I_i$. Since the probability that a given position is not selected in a single selection is $1-1/n$, we have $\E[I_i]=\left(1-1/n\right)^{3k}$. Hence, $\mu_n=\E[Y]=n\left(1-1/n\right)^{3k}$. Note that $\log\left(1-1/n\right)=-1/n+O\left(1/n^2\right)$, so $\log\mu_n=3c+o(1)$. Consequently, as $n\to\infty$, $\mu_n\to e^{3c}$. Moreover, for $i\neq j$, $\E[I_iI_j]=\left(1-2/n\right)^{3k}\leq\left(1-1/n\right)^{6k}=\E[I_i]\E[I_j]$
and hence $\operatorname{Cov}(I_i,I_j)\leq0$. Thus, 
\[
\operatorname{Var}(Y)=\sum_{i=1}^n\operatorname{Var}(I_i)+2\sum_{1\leq i<j\leq n}\operatorname{Cov}(I_i,I_j)
\leq\sum_{i=1}^n\operatorname{Var}(I_i)
\leq\sum_{i=1}^n\E[I_i]
=\mu_n.
\]

Let $M=\left\lfloor\mu_n/2\right\rfloor$ and define $A_M=\{\sigma\in S_n:\operatorname{Fix}(\sigma)\geq M\}$, where $\mathrm{Fix}(\sigma)=\#\{i:\sigma(i)=i\}$ denotes the number of fixed points of $\sigma\in S_n$. Since every position that has never been selected remains fixed, we have
\[
T_3^{*k}(A_M)\geq\Prob(Y\geq M).
\]
Since $M\leq\mu_n/2$, by Chebyshev's inequality, 
\[
T_3^{*k}(A_M^c)
\leq\Prob(Y<M)
\leq\Prob\left(|Y-\mu_n|\geq\frac{\mu_n}{2}\right)
\leq\frac{4\operatorname{Var}(Y)}{\mu_n^2}
\leq\frac4{\mu_n}.
\]
Thus,
\[
T_3^{*k}(A_M)\geq1-\frac4{\mu_n}.
\]
On the other hand, as shown in Section~4 of Chapter~IV of Feller (1968)~\cite{book1968}, we have
\[
U_n(\operatorname{Fix}(\sigma)=r)=\frac{1}{r!}\sum_{j=0}^{n-r}\frac{(-1)^j}{j!}
\]
and hence
\[
U_n(A_M)=\sum_{r=M}^n\frac{1}{r!}\sum_{j=0}^{n-r}\frac{(-1)^j}{j!}\leq\sum_{r=M}^n\frac{1}{r!}\leq\frac{2}{M!}\leq\frac{8}{\mu_n}. 
\]
Therefore,
\[
\dtv(T_3^{*k}, U_n)\geq T_3^{*k}(A_M)-U_n(A_M)\geq1-\frac{12}{\mu_n}. 
\]

Combining these estimates, we conclude that there exists a constant $b>0$, independent of $c$, such that, for all sufficiently large $n$,
\[
\dtv(T_3^{*k}, U_n) \geq1-be^{-3c}.
\]
This completes the proof. 
\end{proof}

\begin{rem*}
Diaconis and Shahshahani (1981)~\cite{paper1981} proved in Theorem 2 that
\[
\dtv(T_2^{*k},U_n)\geq \frac1e-\exp(-e^{-2c})+o(1), \qquad c=\frac{k-\frac12 n\log n}{n}.
\]
However, this lower bound does not by itself establish cutoff for the random transposition shuffle. In their proof, they considered the set of permutations with at least one fixed point. As in the proof of Theorem~\ref{thm:lower-bound}, we can strengthen this bound by considering permutations with at least $M$ fixed points, where $M$ is chosen to be approximately half the expected number of positions that have never been selected during the first $k$ steps. This yields
\[
\dtv(T_2^{*k},U_n) \geq 1-be^{-2c}, \qquad k=\left\lfloor\frac12 n\log n-cn\right\rfloor, \qquad c>0,
\]
for some constant $b>0$, independent of $c$, and all sufficiently large $n$. Together with the upper bound in Theorem 1 of Diaconis and Shahshahani (1981)~\cite{paper1981}, this establishes cutoff at time $\frac12 n\log n$ with a window of size $O(n)$.
\end{rem*}
\vskip \baselineskip

\section{Limit profile (Proof of Theorem \ref{thm:limit-profile})}

In this section, we prove Theorem~\ref{thm:limit-profile}. We first establish some auxiliary lemmas needed for the proof. The following approximation follows from Lemma~\ref{lem2020}. 

\begin{lem}\label{lem5.1}
For every integer $M\geq1$ and sufficiently large $n$, let $S_M(n)$ be a subset of $\widehat{S_n}^*$ defined by
\[
S_M(n)=\left\{\lambda\in\widehat{S_n}^*:\lambda_1\geq n-M\right\}. 
\]
Then
\[
\left|d_1(T_3^{*k}, U_n)-\frac{1}{n!}\sum_{\sigma\in S_n}\left|\sum_{\lambda\in S_M(n)}d_\lambda s_\lambda^k\operatorname{ch}^\lambda(\sigma)\right|\right|\leq\sum_{\lambda\in\widehat{S_n}:\lambda_1<n-M}d_\lambda|s_\lambda|^k. 
\]
\end{lem}

\begin{lem}\label{lem5.2}
For every $\epsilon>0$, there exist $M=M(c,\epsilon)\geq1$ and $n_0=n_0(M)$ such that
\[
\sum_{\lambda\in\widehat{S_n}:\lambda_1\leq n-M}d_\lambda|s_\lambda|^k\leq\epsilon \qquad \text{for all $n\geq n_0$}. 
\]
\end{lem}

\begin{proof}
We first show that $\sum_{\lambda_1\leq n-1}d_\lambda|s_\lambda|^k$ is bounded uniformly in $n$. We use the same zones as in the proof of Theorem~\ref{thm:upper-bound}.

\noindent\textbf{Zone 1}

As in the proof of the upper bound, we have $|s_\lambda| \leq 1/9$ for all sufficiently large $n$. By (\ref{d1-2}) and the definition (\ref{k-2}) of $k$, 
\begin{align*}
\sum_{\lambda:\text{Zone 1}}d_\lambda|s_\lambda|^k
&\leq p(n)\sqrt{n!}\left(\frac19\right)^k\\
&\leq\exp \left(\left(\frac12-\frac23\log 3\right)n\log n+O(n)\right)\\
&\xrightarrow[n\to\infty]{}0.
\end{align*}

\noindent\textbf{Zone 2}

Recall that $|s_\lambda|\leq\dfrac13+o(1)$ and $d_\lambda \leq 2^n n^{n/3}$. Therefore,
\begin{align*}
\sum_{\lambda:\text{Zone 2}}d_\lambda|s_\lambda|^k
&\leq p(n)2^n n^{\frac{n}{3}}\left(\frac13+o(1)\right)^k\\
&\leq\exp \left(\frac13\left(1+\log\left(\frac13+o(1)\right)\right)n\log n+O(n)\right)\\
&\xrightarrow[n\to\infty]{}0.
\end{align*}

\noindent\textbf{Zone 3-(1)}

Write $\lambda_1=n-j$. As in the proof of the upper bound, 
\[
|s_\lambda|\leq 1-3\left(\frac{j}{n}-\frac{j^2}{n^2}+\frac{1}{n^2}\right)
\leq\exp\left(-3\left(\frac{j}{n}-\frac{j^2}{n^2}+\frac{1}{n^2}\right)\right).
\]
Thus,
\begin{align*}
\sum_{\lambda:\text{Zone 3-(1)}}d_\lambda|s_\lambda|^k
&=\sum_{j=1}^{\left\lfloor\frac{n}{3}\right\rfloor}\sum_{\lambda_1=n-j}d_\lambda|s_\lambda|^k\\
&\leq \sum_{j=1}^{\left\lfloor\frac{n}{1000}\right\rfloor}\sum_{\lambda_1=n-j}d_\lambda|s_\lambda|^k+\sum_{j=\left\lfloor\frac{n}{1000}\right\rfloor+1}^{\left\lfloor\frac{n}{3}\right\rfloor}\sum_{\lambda_1=n-j}d_\lambda|s_\lambda|^k\\
&=A_1+A_2.
\end{align*}
By (\ref{d2-1}) and the Cauchy-Schwarz inequality, 
\begin{align}\label{d*}
\sum_{\lambda_1=n-j}d_\lambda
&\leq\sum_{\lambda^{*}\vdash j}\binom{n}{j}d_{\lambda^{*}}
\leq \binom{n}{j}\sqrt{\left(\sum_{\lambda^{*}\vdash j}1^2\right)\left(\sum_{\lambda^{*}\vdash j}d_{\lambda^{*}}^2\right)}\notag\\
&\leq \frac{n^j}{j!}\sqrt{2^jj!}
= \frac{n^j2^{\frac{j}{2}}}{\sqrt{j!}}.
\end{align}
Using (\ref{d*}), we bound $A_1$ as follows:
\begin{align*}
A_1&\leq\sum_{j=1}^{\left\lfloor\frac{n}{1000}\right\rfloor}\frac{n^j2^{\frac{j}{2}}}{\sqrt{j!}}\exp\left(-3k\left(\frac{j}{n}-\frac{j^2}{n^2}+\frac{1}{n^2}\right)\right)\\
&\leq\sum_{j=1}^{\left\lfloor\frac{n}{1000}\right\rfloor}\frac{2^{\frac{j}{2}}}{\sqrt{j!}}\exp\left(j\log n-3\left(\frac{1}{3}n\log n+cn-1\right)\left(\frac{j}{n}-\frac{j^2}{n^2}+\frac{1}{n^2}\right)\right)\\
&\leq\sum_{j=1}^{\left\lfloor\frac{n}{1000}\right\rfloor}\frac{2^{\frac{j}{2}}}{\sqrt{j!}}\exp\left(-3jc+\frac{(j^2-1)(\log n+3c)+3j}{n}\right).
\end{align*}
Let $a_j(n)$ denote the $j$-th summand in the expression above and set
\[
b_j(n)=\frac{a_{j+1}(n)}{a_j(n)}=\frac{\sqrt{2}e^{-3c}}{\sqrt{j+1}}\exp\left(\frac{(2j+1)(\log n+3c)+3}{n}\right). 
\]
The sequence $b_j(n)$ first decreases and then increases, with a turning point near $j=\dfrac{n}{4(\log n+3c)}$. Moreover, evaluating the displayed expression at $j=\lfloor n/1000\rfloor$ gives
\begin{align*}
b_{\left\lfloor\frac{n}{1000}\right\rfloor}(n)&\leq\frac{\sqrt{2}e^{-3c}}{\sqrt{\frac{n}{1000}}}\exp\left(\frac{\left(\frac{2n}{1000}+1\right)(\log n+3c)+3}{n}\right)\\
&=\sqrt{2000}\exp\left(-3c+\frac{6c}{1000}+\frac{\log n+3c+3}{n}\right)n^{\frac{2}{1000}-\frac12}\\
&\xrightarrow[n\to\infty]{}0.
\end{align*}
Choose $M$ such that
\begin{equation}\label{lem5.2-2}
\frac{\sqrt{2}e^{-3c}}{\sqrt{M+1}}\leq \frac14
\end{equation}
and take $n$ sufficiently large so that
\[
\exp\left(\frac{(2M+1)(\log n+3c)+3}{n}\right)\leq 2. 
\]
Then
\[
A_1\leq \sum_{j=1}^Ma_j(n)+a_M(n)\sum_{i=1}^{\infty}\frac{1}{2^i}
\xrightarrow[n\to\infty]{}\sum_{j=1}^M\frac{2^{\frac{j}{2}}}{\sqrt{j!}}e^{-3jc}+\frac{2^{\frac{M}{2}}}{\sqrt{M!}}e^{-3Mc}.
\]
For the remaining range, we use $j!\geq(j/e)^j$ and (\ref{d*}) to obtain
\[
\sum_{\lambda_1=n-j}d_\lambda\leq\frac{n^j(2e)^{\frac{j}{2}}}{j^{\frac{j}{2}}}.
\]
Since $1\leq j\leq n/3$, $|s_\lambda|\leq1-2j/n\leq e^{-2j/n}$. Therefore, for $n/1000<j\leq n/3$,
\begin{align*}
\sum_{\lambda_1=n-j}d_\lambda|s_\lambda|^k
&\leq\frac{n^j(2e)^{\frac{j}{2}}}{j^{\frac{j}{2}}}e^{-\frac{2jk}{n}}\\
&\leq\exp\left(\frac{j}{3}\log n-\frac{j}{2}\log j+\frac{j}{2}\log2e-2cj+\frac{2j}{n}\right)\\
&\leq\exp\left(-\frac{j}{6}\log n+\frac{j}{2}\log2000e+2|c|j+\frac23\right)\\
&\leq\exp\left(-\frac{n}{6000}\log n+O(n)\right). 
\end{align*}
Hence,
\[
A_2\leq\frac{n}{3}\exp\left(-\frac{n}{6000}\log n+O(n)\right)\xrightarrow[n\to\infty]{}0.
\]

\noindent\textbf{Zone 3-(2)}

The argument used for Zone 3-(2) in the proof of the upper bound gives
\[
\sum_{\lambda:\text{Zone 3-(2)}}d_{\lambda}|s_\lambda|^k
\leq \left(\sum_{\lambda:\text{Zone 3-(1)}}d_\lambda|s_\lambda|^k\right)\left(1-\frac{2(n-3)}{n^2}\right)^k+\sum_{\lambda:\text{Zone 3-(1)}}d_\lambda\left(\frac6n\right)^k.
\]
For the first term,
\begin{align*}
&\left(\sum_{\lambda:\text{Zone 3-(1)}}d_\lambda|s_\lambda|^k\right)\left(1-\frac{2(n-3)}{n^2}\right)^k\\
&\leq\left(\sum_{\lambda:\text{Zone 3-(1)}}d_\lambda|s_\lambda|^k\right)\exp\left(-\frac{2}{3}\log n-2c+O\left(\frac{\log n}{n}\right)\right)\\
&\xrightarrow[n\to\infty]{}0. 
\end{align*}
For the second term, 
\[
\sum_{\lambda:\text{Zone 3-(1)}}d_\lambda\left(\frac6n\right)^k
\leq\exp\left(-\frac13n(\log n)^2+O(n\log n)\right)
\xrightarrow[n\to\infty]{}0.
\]
Therefore,
\[
\sum_{\lambda:\text{Zone 3-(2)}}d_{\lambda}|s_\lambda|^k\xrightarrow[n\to\infty]{}0.
\]

\noindent\textbf{Zone 3-(3)}

Since $s_{(1^n)}=1-\dfrac6n+\dfrac{6}{n^2}$, we have
\[
d_{(1^n)}|s_{(1^n)}|^k
\leq\exp\left(-2\log n-6c+O\left(\frac{\log n}{n}\right)\right)
\xrightarrow[n\to\infty]{}0.
\]
\vskip \baselineskip

Combining the estimates above, for any $\epsilon>0$, we choose $M=M(c,\epsilon)\geq1$ such that both (\ref{lem5.2-2}) and 
\[
4\cdot\frac{2^{\frac{M}{2}}}{\sqrt{M!}}e^{-3Mc}<\frac{\epsilon}{2}
\]
hold. Moreover, we take $n$ sufficintly large so that
\[
a_M(n)\leq2\lim_{n\to\infty}a_M(n)=2\cdot\frac{2^{\frac{M}{2}}}{\sqrt{M!}}e^{-3Mc}
\]
as well. Then 
\begin{align*}
&\sum_{\lambda_1\leq n-M}d_\lambda|s_\lambda|^k\\
&\leq \left(\sum_{\lambda:\text{Zone 1}}+\sum_{\lambda:\text{Zone 2}}+\sum_{\lambda:\text{Zone 3-(2)}}+\sum_{\lambda:\text{Zone 3-(3)}}\right)d_\lambda|s_\lambda|^k+\sum_{j=M}^{\left\lfloor\frac{n}{1000}\right\rfloor}a_j(n)+A_2\\
&\leq \sum_{j=M}^{\left\lfloor\frac{n}{1000}\right\rfloor}a_j(n)+\frac{\epsilon}{2}\\
&\leq a_M(n)\sum_{i=0}^{\infty}\frac{1}{2^i}+\frac{\epsilon}{2}\\
&\leq 4\cdot\frac{2^{\frac{M}{2}}}{\sqrt{M!}}e^{-3Mc}+\frac{\epsilon}{2}\\
&<\epsilon. 
\end{align*}
This completes the proof of Lemma~\ref{lem5.2}.  
\end{proof}
\vskip \baselineskip

To complete the proof of Theorem~\ref{thm:limit-profile}, we use the following four results, adapted from Lemmas 4.2, 4.4, and 4.6 and Proposition 4.5 of Teyssier (2020) \cite{paper2020}.

\begin{lem}[{cf.~\cite[Lemma 4.2]{paper2020}}]\label{lem5.3}
Let $l$ be a fixed positive integer. Then, as $n\to\infty$,
\[
\frac{1}{n!}\sum_{\sigma\in S_n}\left|\sum_{j=1}^l\sum_{\lambda_1=n-j}d_\lambda s_\lambda^k\ch^\lambda(\sigma)\right|=\frac{1}{n!}\sum_{\sigma\in S_n}\left|\sum_{j=1}^le^{-3jc}\Phi_j(\operatorname{Fix}(\sigma))\right|+o(1),
\]
where for every nonnegative integer $r$, 
\[
\Phi_j(r)=\sum_{i=0}^j\binom{r}{j-i}\frac{(-1)^i}{i!}, 
\]
with the convention that $\dbinom{a}{b}=0$ whenever $a<b$.
\end{lem}

\begin{lem}[{cf.~\cite[Lemma 4.4]{paper2020}}]\label{lem5.4}
Let $\epsilon>0$. There exists $M_0=M_0(\epsilon,c)$ such that, for every $M\geq M_0$ and every positive integer $n$,
\[
\left|\frac{1}{n!}\sum_{\sigma\in S_n}\left|\sum_{j=1}^Me^{-3jc}\Phi_j(\operatorname{Fix}(\sigma))\right|-\frac{1}{n!}\sum_{\sigma\in S_n}\left|\sum_{j=1}^\infty e^{-3jc}\Phi_j(\operatorname{Fix}(\sigma))\right|\right|\leq\epsilon.
\]
\end{lem}

Define the function $f_c$ by
\[
f_c(x)=e^{-e^{-3c}}(1+e^{-3c})^x-1, \qquad x\in\R.
\]

\begin{prop}[{cf.~\cite[Proposition 4.5]{paper2020}}]\label{prop5.5}
For a positive integer $N$,
\[
\sum_{j=1}^\infty e^{-3jc}\Phi_j(N)=f_c(N).
\]
\end{prop}

\begin{lem}[{cf.~\cite[Lemma 4.6]{paper2020}}]\label{lem5.6}
As $n\to\infty$,
\[
\frac{1}{n!}\sum_{\sigma\in S_n}|f_c(\operatorname{Fix}(\sigma))|\to\E[|f_c(\operatorname{Poiss}(1))|]. 
\]
\end{lem}

The proof of Theorem~\ref{thm:limit-profile} proceeds along the same lines as that of {\cite[Theorem 1.1]{paper2020}}.

\begin{proof}[\normalfont\bfseries Proof of Theorem \ref{thm:limit-profile}]

Fix $\epsilon>0$. Choose $M$ sufficiently large and then choose $n_0$ so that the estimates below hold for all $n\geq n_0$. By Lemmas~\ref{lem5.1} and~\ref{lem5.2}, 
\[
\left|d_1(T_3^{*k}, U_n)-\frac{1}{n!}\sum_{\sigma\in S_n}\left|\sum_{\lambda\in S_M(n)}d_\lambda s_\lambda^k\ch^\lambda(\sigma)\right|\right|\leq\epsilon.
\]
By Lemma~\ref{lem5.3},
\[
\left|\frac{1}{n!}\sum_{\sigma\in S_n}\left|\sum_{\lambda\in S_M(n)}d_\lambda s_\lambda^k\ch^\lambda(\sigma)\right|-\frac{1}{n!}\sum_{\sigma\in S_n}\left|\sum_{j=1}^M e^{-3jc}\Phi_j(\operatorname{Fix}(\sigma))\right|\right|\leq\epsilon.
\]
By Lemma~\ref{lem5.4},
\[
\left|\frac{1}{n!}\sum_{\sigma\in S_n}\left|\sum_{j=1}^M e^{-3jc}\Phi_j(\operatorname{Fix}(\sigma))\right|-\frac{1}{n!}\sum_{\sigma\in S_n}\left|\sum_{j=1}^\infty e^{-3jc}\Phi_j(\operatorname{Fix}(\sigma))\right|\right|\leq\epsilon.
\]
By Proposition~\ref{prop5.5} and Lemma~\ref{lem5.6},
\begin{align*}
&\left|\frac{1}{n!}\sum_{\sigma\in S_n}\left|\sum_{j=1}^\infty e^{-3jc}\Phi_j(\operatorname{Fix}(\sigma))\right|-\E[|f_c(\operatorname{Poiss}(1))|]\right|\\
&=\left|\frac{1}{n!}\sum_{\sigma\in S_n}\left|f_c(\operatorname{Fix}(\sigma))\right|-\E[|f_c(\operatorname{Poiss}(1))|]\right|\\
&\leq\epsilon.
\end{align*}
Thus, by the triangle inequality,
\[
|d_1(T_3^{*k}, U_n)-\E[|f_c(\operatorname{Poiss}(1))|]|\leq4\epsilon.
\]
Consequently, for every $c\in\R$,
\[
d_1(T_3^{*k}, U_n)\xrightarrow[n\to\infty]{}\E[|f_c(\operatorname{Poiss}(1))|].
\]
Finally,
\begin{align*}
\E[|f_c(\operatorname{Poiss}(1))|]&=\sum_{r=0}^\infty\frac{e^{-1}}{r!}|e^{-e^{-3c}}(1+e^{-3c})^r-1|\\
&=\sum_{r=0}^\infty\left|\frac{(e^{1+e^{-3c}})^{-1}}{r!}(1+e^{-3c})^r-\frac{e^{-1}}{r!}1^r\right|\\
&=d_1(\operatorname{Poiss}(1+e^{-3c}), \operatorname{Poiss}(1))
\end{align*}
and the desired result follows.
\end{proof}
\vskip \baselineskip

\section{Biased random transposition-3-cycle shuffle}\label{Biased random transposition-3-cycle shuffle}

We now consider a biased version of the random transposition-3-cycle shuffle, in which the relative frequencies of transpositions and 3-cycles are controlled by a parameter $p\in(0,1)$. At each step, the identity permutation is chosen with probability $1/n^2$. Conditionally on a nonidentity step, a transposition is chosen with probability $1-p$ and a 3-cycle with probability $p$. More precisely, the probability distribution is defined as follows.

\begin{dfn}
Let $p\in(0,1)$. Define a probability $T_p$ on the symmetric group $S_n$ by
\[
T_p(\sigma)=
\begin{cases}
\dfrac{1}{n^2},
   & \sigma=\mathrm{id},\\[1ex]
(1-p)\cdot \dfrac{2(n+1)}{n^3},
   & \sigma=(i\ j),\\[1ex]
p\cdot \dfrac{3(n+1)}{n^3(n-2)},
   & \sigma=(i\ j\ k),\\[1ex]
0,
   & \text{otherwise}. 
\end{cases}
\]
The random walk on $S_n$ driven by $T_p$ is called the biased random transposition-3-cycle shuffle.
\end{dfn}

Applying Lemma~\ref{lem1} to $T_p$, we obtain
\begin{equation}\label{s_p_lambda}
s_\lambda=\frac{1}{n^2}+(1-p)\left(1-\frac{1}{n^2}\right)r_2(\lambda)+p\left(1-\frac{1}{n^2}\right)r_3(\lambda).
\end{equation}

\subsection{Upper bound}

\begin{thm}\label{thm:biased1}
For $p\in(0,1)$, let
\[
k=\left\lfloor \frac{1}{2+p}n\log n+cn\right\rfloor, \qquad c>0. 
\]
Then there exists a constant $b>0$, independent of $c$, such that, for all sufficiently large $n$,
\[
\dtv(T_p^{*k}, U_n) \leq be^{-(2+p)c}. 
\]
\end{thm}

\begin{proof}
Applying Lemma~\ref{lem1981} to $T_p$, we obtain
\[
d_1(T_p^{*k}, U_n)^2 \leq \sum_{\lambda\in\widehat{S_n}^*} d_{\lambda}^2 s_{\lambda}^{2k}.
\]

We estimate the sum on the right-hand side by considering the following zones: 

\begin{align*}
&\text{Zone 1. }\qquad  1 \leq \lambda_1 \leq \frac{n}{3} \ \text{and} \ 1 \leq m \leq \frac{n}{3}; \\
&\text{Zone 2. }\qquad  \left(\frac{n}{3}<\lambda_1\leq\frac{n}{2} \ \text{and} \ m\leq\frac{n}{2}\right) \ \text{or} \ \left(\lambda_1 \leq \frac{n}{2} \ \text{and} \ \frac{n}{3}<m\leq\frac{n}{2}\right); \\
&\text{Zone 3. }\qquad  \frac{n}{2}<\lambda_1<\frac{2}{3}n \ \text{or} \ \frac{n}{2}<m<\frac{2}{3}n; \\
&\text{Zone 4. }\qquad  (1)\ \frac{2}{3}n \leq \lambda_1 \leq n-1 \qquad (2)\ \frac{2}{3}n \leq m \leq n-1 \qquad (3)\ m=n. \\
\end{align*}

\begin{center}
\begin{tikzpicture}[scale=4, xscale=1.5]
\draw[->] (0,0) -- (1.1,0)
  node[right] {$\lambda_1$};

\draw[->] (0,0) -- (0,1.1)
  node[above] {$m$};

\draw (0.05,0.05) -- (0.05,1);
\draw (1/3,0.05) -- (1/3,1/3);
\draw (1/2,0.05) -- (1/2,1/2);
\draw (2/3,0.05) -- (2/3,1/3+0.05);
\draw (0.05,0.05) -- (1,0.05);
\draw (0.05,1/3) -- (1/3,1/3);
\draw (0.05,1/2) -- (1/2,1/2);
\draw (0.05,2/3) -- (1/3+0.05,2/3);
\draw (0.05,1) -- (1,0.05);

\node[below] at (0.05,0) {$1$};
\node[below] at (1/3,0) {$n/3$};
\node[below] at (1/2,0) {$n/2$};
\node[below] at (2/3,0) {$2n/3$};
\node[below] at (1,0) {$n$};
\node[left] at (0,0.05) {$1$};
\node[left] at (0,1/3) {$n/3$};
\node[left] at (0,1/2) {$n/2$};
\node[left] at (0,2/3) {$2n/3$};
\node[left] at (0,1.) {$n$};

\draw[->] (0.2,1.05) -- (0.05,1);
\node[right] at (0.2,1.05) {4-(3)};

\node at (1/6+0.025,1/6+0.025) {1};
\node at (0.4,0.4) {2};
\node at (1/4,7/12) {3};
\node at (7/12,1/4) {3};
\node at (7/9,1/12+0.025) {4-(1)};
\node at (1/9+0.025,3/4) {4-(2)};
\end{tikzpicture}
\end{center}

\noindent\textbf{Zone 1}

By Lemmas~\ref{r2-1}, \ref{r3-1}, \ref{r3-5} and (\ref{r2-lambda'}), we have
\[
-\frac13\leq r_2(\lambda) \leq \frac13, \qquad
-\frac{3}{2(n-2)} \leq r_3(\lambda) \leq \frac19. 
\]
(\ref{s_p_lambda}) implies $|s_\lambda|\leq \dfrac13-\dfrac29p+o(1)$. Therefore, 
\begin{align*}
\sum_{\lambda:\text{Zone 1}}d_\lambda^2s_\lambda^{2k}
&\leq \left(\frac13-\frac29p+o(1)\right)^{2k}n^n\\
&\leq\exp \left(\left(1+\frac{2}{2+p}\log\left(\frac13-\frac29p+o(1)\right)\right)n\log n+O(n)\right)\\
&\xrightarrow[n\to\infty]{}0.
\end{align*}

\noindent\textbf{Zone 2}

By Lemmas~\ref{r2-2}, \ref{r3-2}, \ref{r3-5} and (\ref{r2-lambda'}), we have
\[
-\frac12\leq r_2(\lambda) \leq \frac12, \qquad
0 \leq r_3(\lambda) \leq \frac14. 
\]
By (\ref{s_p_lambda}), $|s_\lambda|\leq \dfrac12-\dfrac{p}{4}+o(1)$. (\ref{d2-2}) implies $d_\lambda \leq 2^n n^{n/3}$, so
\begin{align*}
\sum_{\lambda:\text{Zone 2}}d_\lambda^2s_\lambda^{2k}
&\leq p(n)4^nn^{\frac23n}\left(\frac12-\frac{p}{4}+o(1)\right)^{2k}\\
&\leq\exp\left(\left(\frac23+\frac{2}{2+p}\log\left(\frac12-\frac{p}{4}+o(1)\right)\right)n\log n+O(n)\right)\\
&\xrightarrow[n\to\infty]{}0.
\end{align*}

\noindent\textbf{Zone 3}

By Lemmas~\ref{r2-3}, \ref{r3-3}, \ref{r3-5} and (\ref{r2-lambda'}), we have
\[
-\frac59 \leq r_2(\lambda) \leq \frac59, \qquad 0 \leq r_3(\lambda) \leq \frac13. 
\]
By (\ref{s_p_lambda}), $|s_\lambda|\leq \dfrac59-\dfrac29p+o(1)$. (\ref{d2-2}) implies $d_\lambda \leq 2^n n^{n/4}$, so
\begin{align*}
\sum_{\lambda:\text{Zone 3}}d_\lambda^2s_\lambda^{2k}
&\leq p(n)4^nn^{\frac{n}{2}}\left(\frac59-\frac29p+o(1)\right)^{2k}\\
&\leq\exp \left(\left(\frac12+\frac{2}{2+p}\log\left(\frac59-\frac29p+o(1)\right)\right)n\log n+O(n)\right)\\
&\xrightarrow[n\to\infty]{}0.
\end{align*}

\noindent\textbf{Zone 4-(1)}

Write $\lambda_1=n-j$. By Lemmas~\ref{r2-4}, \ref{r3-4} and (\ref{s_p_lambda}), we obtain
\begin{align*}
s_\lambda
&\leq\frac{1}{n^2}+(1-p)\left(1-\frac{1}{n^2}\right)\left(1-\frac{2j(n-j+1)}{n(n-1)}\right)\\
&\quad+p\left(1-\frac{1}{n^2}\right)\left(1-\frac{3j(n-j-1)}{(n-1)(n-2)}\right)\\
&=1-(2+p)\left(\frac{j}{n}-\frac{j^2}{n^2}+\frac{j}{n^2}\right)\\
&\quad-\frac{2(1-p)}{n}\left(\frac{j}{n}-\frac{j^2}{n^2}+\frac{j}{n^2}\right)-\frac{3pj(n+1-3j)}{n^2(n-2)}\\
&\leq1-(2+p)\left(\frac{j}{n}-\frac{j^2}{n^2}+\frac{j}{n^2}\right)\\
&\leq\exp\left(-(2+p)\left(\frac{j}{n}-\frac{j^2}{n^2}+\frac{j}{n^2}\right)\right).
\end{align*}
Since $r_2(\lambda)\geq0$ and $r_3(\lambda)\geq0$ by Lemmas~\ref{r2-6} and~\ref{r3-5}, respectively, for all sufficiently large $n$, we have $s_\lambda\geq0$. Therefore,
\[
|s_\lambda|\leq\exp\left(-(2+p)\left(\frac{j}{n}-\frac{j^2}{n^2}+\frac{j}{n^2}\right)\right).
\]
By (\ref{d3}), 
\begin{align*}
\sum_{\lambda:\text{Zone 4-(1)}}d_\lambda^2s_\lambda^{2k}
&=\sum_{j=1}^{\left\lfloor\frac{n}{3}\right\rfloor}\sum_{\lambda_1=n-j}d_\lambda^2s_\lambda^{2k}\\
&\leq\sum_{j=1}^{\left\lfloor\frac{n}{3}\right\rfloor}p(j)\frac{n^{2j}}{j!}\exp\left(-2k(2+p)\left(\frac{j}{n}-\frac{j^2}{n^2}+\frac{j}{n^2}\right)\right)\\
&\leq\sum_{j=1}^{\left\lfloor\frac{n}{3}\right\rfloor}\frac{p(j)}{j!}\exp\left(\frac{2\log n}{n}j(j-1)-2(2+p)(cn-1)\left(\frac{j}{n}-\frac{j^2}{n^2}+\frac{j}{n^2}\right)\right)\\
&\leq e^{-2(2+p)c+\frac{2(2+p)}{3}}\sum_{j=1}^{\left\lfloor\frac{n}{3}\right\rfloor}\frac{p(j)}{j!}\exp\left(\frac{2\log n}{n}j(j-1)\right).
\end{align*}
By (\ref{S(n)}), 
\[
\sup_{n\geq1}\sum_{j=1}^{\left\lfloor\frac{n}{3}\right\rfloor}\frac{p(j)}{j!}\exp\left(\frac{2\log n}{n}j(j-1)\right)\leq\sup_{n\geq1}S(n)<\infty. 
\]
Hence, there exists a constant $B>0$, independent of $c$, such that, for all sufficiently large $n$, 
\[
\sum_{\lambda:\text{Zone 4-(1)}}d_\lambda^2s_\lambda^{2k} \leq Be^{-2(2+p)c}.
\]

\noindent\textbf{Zone 4-(2)}

Let $\lambda'$ be the conjugate partition of $\lambda$. By (\ref{r2-lambda'}) and (\ref{r3-lambda'}), we have
\begin{align*}
s_\lambda
&=\frac{1}{n^2}-(1-p)\left(1-\frac{1}{n^2}\right)r_2(\lambda')+p\left(1-\frac{1}{n^2}\right)r_3(\lambda')\\
&=s_{\lambda'}-2(1-p)\left(1-\frac{1}{n^2}\right)r_2(\lambda').
\end{align*}
Since $2n/3\leq\lambda'_1\leq n-1$, Lemma~\ref{r2-5} implies that
\begin{align*}
s_\lambda
&\leq s_{\lambda'}-2(1-p)\left(1-\frac{1}{n^2}\right)\frac{n-3}{3(n-1)}\\
&=s_{\lambda'}-2(1-p)\frac{(n+1)(n-3)}{3n^2}.
\end{align*}
Set $\delta_n=2(1-p)\dfrac{(n+1)(n-3)}{3n^2}$. If $s_\lambda\geq0$, then $s_{\lambda'}\geq\delta_n>0$ and $s_{\lambda'}\leq1$. Hence,
\begin{equation}\label{thm6.2-1}
|s_\lambda|=s_\lambda\leq s_{\lambda'}-\delta_n\leq(1-\delta_n)|s_{\lambda'}|.
\end{equation}
On the other hand, if $s_\lambda<0$, let $\lambda'_1=n-j$. By Lemmas~\ref{r2-4} and~\ref{r3-5},
\begin{equation}\label{thm6.2-2}
|s_\lambda|=-s_\lambda\leq(1-p)\left(1-\frac{2j}{n}+\frac{2j^2}{n^2}\right).
\end{equation}
By (\ref{thm6.2-1}),
\begin{align*}
\sum_{\lambda:\text{Zone 4-(2)}}d_\lambda^2s_\lambda^{2k}
&\leq\sum_{\substack{\lambda:\text{Zone 4-(2)}\\s_\lambda\geq0}}d_\lambda^2(1-\delta_n)^{2k}s_{\lambda'}^{2k}+\sum_{\substack{\lambda:\text{Zone 4-(2)}\\s_\lambda<0}}d_\lambda^2s_\lambda^{2k}\\
&\leq\left(\sum_{\lambda:\text{Zone 4-(1)}}d_\lambda^2s_\lambda^{2k}\right)(1-\delta_n)^{2k}+\sum_{\substack{\lambda:\text{Zone 4-(2)}\\s_\lambda<0}}d_\lambda^2s_\lambda^{2k}.
\end{align*}
For the first term,
\begin{align*}
&\left(\sum_{\lambda:\text{Zone 4-(1)}}d_\lambda^2s_\lambda^{2k}\right)(1-\delta_n)^{2k}\\
&\leq\left(\sum_{\lambda:\text{Zone 4-(1)}}d_\lambda^2s_\lambda^{2k}\right)\exp\left(-4(1-p)\frac{(n+1)(n-3)}{3n^2}\left(\frac{1}{2+p}n\log n+cn-1\right)\right)\\
&=\left(\sum_{\lambda:\text{Zone 4-(1)}}d_\lambda^2s_\lambda^{2k}\right)\exp\left(-\frac{4(1-p)}{3(2+p)}n\log n+O(n)\right)\\
&\xrightarrow[n\to\infty]{}0.
\end{align*}
For the second term, by (\ref{lambda'}), (\ref{d3}) and (\ref{thm6.2-2}), 
\[
\sum_{\substack{\lambda:\text{Zone 4-(2)}\\s_\lambda<0}}d_\lambda^2s_\lambda^{2k}
\leq\sum_{j=1}^{\left\lfloor\frac{n}{3}\right\rfloor}p(j)\frac{n^{2j}}{j!}\left((1-p)\left(1-\frac{2j}{n}+\frac{2j^2}{n^2}\right)\right)^{2k}.
\]
By (\ref{p(n)-2}) and $j!\geq(j/e)^j$, 
\begin{align*}
\log\left(p(j)\frac{n^{2j}}{j!}\right)
&\leq2j\log n-j\log j+j+\pi\left(\frac23j\right)^{\frac12}\\
&\leq j\log n+O(n).
\end{align*}
Moreover,
\[
\log\left(1-\frac{2j}{n}+\frac{2j^2}{n^2}\right)\leq-\frac{3j}{2n}. 
\]
Hence,
\begin{align*}
&\log\left(p(j)\frac{n^{2j}}{j!}\left((1-p)\left(1-\frac{2j}{n}+\frac{2j^2}{n^2}\right)\right)^{2k}\right)\\
&\leq j\log n+2\left(\frac{1}{2+p}n\log n+cn-1\right)\left(\log(1-p)-\frac{3j}{2n}\right)+O(n)\\
&=\frac{2}{2+p}\log(1-p)\cdot n\log n-\frac{1-p}{2+p}j\log n+O(n)\\
&\leq\frac{2}{2+p}\log(1-p)\cdot n\log n+O(n).
\end{align*}
Therefore,
\[
\sum_{\substack{\lambda:\text{Zone 4-(2)}\\s_\lambda<0}}d_\lambda^2s_\lambda^{2k}
\leq\frac n3\exp\left(\frac{2}{2+p}\log(1-p)\cdot n\log n+O(n)\right)
\xrightarrow[n\to\infty]{}0.
\]
Consequently,
\[
\sum_{\lambda:\text{Zone 4-(2)}}d_\lambda^2s_\lambda^{2k}\xrightarrow[n\to\infty]{}0.
\]

\noindent\textbf{Zone 4-(3)}

Since $r_2((1^n))=-1$ and $r_3((1^n))=1$, we have
\[
s_{(1^n)}=1-2(1-p)\left(1-\frac{1}{n^2}\right)\xrightarrow[n\to\infty]{}2p-1.
\]
Since $p\in(0,1)$, we have $|2p-1|<1$. Thus, there exists a constant $q_p\in(0,1)$ such that $|s_{(1^n)}|\leq q_p$ for all sufficiently large $n$. Therefore, 
\[
d_{(1^n)}^2s_{(1^n)}^{2k}\leq q_p^{2k}\xrightarrow[n\to\infty]{}0.
\]
\vskip \baselineskip

Combining these estimates, we conclude that there exists a constant $B'>0$, independent of $c$, such that, for all sufficiently large $n$, 
\[
d_1(T_p^{*k}, U_n)^2 \leq \sum_{\lambda:\text{Zone1-4}}d_{\lambda}^2s_\lambda^{2k} \leq B'e^{-2(2+p)c}.
\]
The desired result follows from this.
\end{proof}
\vskip \baselineskip

\subsection{Lower bound}

\begin{thm}\label{thm:biased-lower-bound}
For $p\in(0,1)$, let
\[
k=\left\lfloor\frac{1}{2+p}n\log n-cn\right\rfloor, \qquad c>0. 
\]
Then there exists a constant $b>0$, independent of $c$, such that, for all sufficiently large $n$,
\[
\dtv(T_p^{*k},U_n)\geq1-be^{-(2+p)c}. 
\]
\end{thm}

\begin{proof}
Let $Y$ denote the number of positions that have never been involved in a move during the first $k$ steps. For each $i\in\{1,\ldots,n\}$, let $I_i$ be the indicator that position $i$ has never been involved in a move. Then $Y=\sum_{i=1}^n I_i$. At each step, the identity permutation is chosen with probability $1/n^2$, while a transposition and a 3-cycle are chosen with probabilities $(1-p)\cdot(n^2-1)/n^2$ and $p\cdot(n^2-1)/n^2$, respectively. Conditional on choosing a transposition, a given position $i$ is involved with probability $2/n$; conditional on choosing a 3-cycle, the corresponding probability is $3/n$. Hence, the probability that position $i$ is not involved in a move in a single step is 
\[
q_{n,p}=1-\dfrac{(2+p)(n^2-1)}{n^3}. 
\]
Hence, $\E[I_i]=q_{n,p}^{\,k}$, so $\mu_n=\E[Y]=nq_{n,p}^{\,k}$. Therefore, $\log\mu_n=(2+p)c+o(1)$. Consequently, as $n\to\infty$, $\mu_n\to e^{(2+p)c}$.

Fix distinct positions $i$ and $j$. Let $q_{n,p}^{(2)}$ denote the probability that neither position is involved in a move during a single step. The probability that a transposition involves neither $i$ nor $j$ is $\dbinom{n-2}{2}\Bigm/\dbinom{n}{2}$. The probability that a 3-cycle involves neither $i$ nor $j$ is $\dbinom{n-2}{3}\Bigm/\dbinom n3$. Hence,
\begin{align*}
q_{n,p}^{(2)}&=\frac1{n^2}+\frac{n^2-1}{n^2}\left((1-p)\frac{(n-2)(n-3)}{n(n-1)}+p\frac{(n-3)(n-4)}{n(n-1)}\right). 
\end{align*}
Hence,
\[
q_{n,p}^{\,2}-q_{n,p}^{(2)}=\frac{n+1}{n^6}\left((p^2+2)n^3-(p+2)^2(n^2+n-1)\right).
\]
Since $p\in(0,1)$, we have $p^2+2\geq2$ and $(p+2)^2\leq9$. Thus, for all sufficiently large $n$, $q_{n,p}^{(2)}\leq q_{n,p}^{\,2}$. Hence, $\E[I_iI_j]=\left(q_{n,p}^{(2)}\right)^k\leq q_{n,p}^{\,2k}=\E[I_i]\E[I_j]$. Therefore, $\operatorname{Cov}(I_i,I_j)\leq0$, so $\operatorname{Var}(Y)\leq\sum_{i=1}^{n}\operatorname{Var}(I_i)\leq\mu_n$. 

Let $M=\left\lfloor\mu_n/2\right\rfloor$ and define $A_M=\{\sigma\in S_n:\operatorname{Fix}(\sigma)\geq M\}$. Then, as in the proof of Theorem~\ref{thm:lower-bound}, $T_p^{*k}(A_M)\geq1-4/{\mu_n}$. On the other hand, as shown in the proof of Theorem~\ref{thm:lower-bound}, $U_n(A_M)\leq8/{\mu_n}$. Hence,
\[
\dtv(T_p^{*k}, U_n)\geq T_p^{*k}(A_M)-U_n(A_M)\geq1-\frac{12}{\mu_n}. 
\]

Combining these estimates, we conclude that there exists a constant $b>0$, independent of $c$, such that, for all sufficiently large $n$,
\[
\dtv(T_p^{*k}, U_n)\geq1-be^{-(2+p)c}.
\]
This completes the proof. 
\end{proof}
\vskip \baselineskip

\subsection{Limit profile}

As in the unbiased case, we determine the limit profile by applying Lemma~\ref{lem2020} to the biased shuffle.

\begin{thm}\label{thm:biased limit-profile}
For $p\in(0,1)$, let
\[
k=\left\lfloor\frac{1}{2+p}n\log n+cn\right\rfloor, \qquad c\in\R. 
\]
Then
\[
\lim_{n\to\infty}\dtv(T_p^{*k}, U_n)=\dtv\left(\operatorname{Poiss}(1+e^{-(2+p)c}),\operatorname{Poiss}(1)\right). 
\]
\end{thm}

The following approximation follows from Lemma~\ref{lem2020}.

\begin{lem}\label{lem6.1}
For every integer $M\geq1$ and sufficiently large $n$, let $S_M(n)$ be a subset of $\widehat{S_n}^*$ defined by
\[
S_M(n)=\left\{\lambda\in\widehat{S_n}^*:\lambda_1\geq n-M\right\}. 
\]
Then
\[
\left|d_1(T_p^{*k}, U_n)-\frac{1}{n!}\sum_{\sigma\in S_n}\left|\sum_{\lambda\in S_M(n)}d_\lambda s_\lambda^k\operatorname{ch}^\lambda(\sigma)\right|\right|\leq\sum_{\lambda\in\widehat{S_n}:\lambda_1<n-M}d_\lambda|s_\lambda|^k. 
\]
\end{lem}

\begin{lem}\label{lem6.2}
For every $\epsilon>0$, there exist $M=M(c,\epsilon)\geq1$ and $n_0=n_0(M)$ such that
\[
\sum_{\lambda_1\leq n-M}d_\lambda|s_\lambda|^k\leq\epsilon \qquad \text{for all $n\geq n_0$}. 
\]
\end{lem}

\begin{proof}
We show that $\sum_{\lambda_1\leq n-1}d_\lambda|s_\lambda|^k$ is uniformly bounded in $n$ using the same zones as in the proof of Theorem~\ref{thm:biased1}. 

\noindent\textbf{Zone 1}

As in the proof of the upper bound, we have $|s_\lambda|\leq\dfrac13-\dfrac29p+o(1)$. Hence,
\begin{align*}
\sum_{\lambda:\text{Zone 1}}d_\lambda|s_\lambda|^k
&\leq p(n)\sqrt{n!}\left(\frac13-\frac29p+o(1)\right)^k\\
&\leq\exp\left(\left(\frac12+\frac{1}{2+p}\log\left(\frac13-\frac29p+o(1)\right)\right)n\log n+O(n)\right)\\
&\xrightarrow[n\to\infty]{}0.
\end{align*}

\noindent\textbf{Zone 2}

Recall that $|s_\lambda|\leq\dfrac12-\dfrac p4+o(1)$ and $d_\lambda\leq2^n n^{n/3}$. Therefore, 
\begin{align*}
\sum_{\lambda:\text{Zone 2}}d_\lambda|s_\lambda|^k
&\leq p(n)2^n n^{\frac n3}\left(\frac12-\frac p4+o(1)\right)^k\\
&\leq\exp\left(\left(\frac13+\frac{1}{2+p}\log\left(\frac12-\frac p4+o(1)\right)\right)n\log n+O(n)\right)\\
&\xrightarrow[n\to\infty]{}0.
\end{align*}

\noindent\textbf{Zone 3}

Recall that $|s_\lambda|\leq\dfrac59-\dfrac29p+o(1)$ and $d_\lambda\leq2^n n^{n/4}$. Therefore, 
\begin{align*}
\sum_{\lambda:\text{Zone 3}}d_\lambda|s_\lambda|^k
&\leq p(n)2^n n^{\frac n4}\left(\frac59-\frac29p+o(1)\right)^k\\
&\leq\exp\left(\left(\frac14+\frac{1}{2+p}\log\left(\frac59-\frac29p+o(1)\right)\right)n\log n+O(n)\right)\\
&\xrightarrow[n\to\infty]{}0.
\end{align*}

\noindent\textbf{Zone 4-(1)}

Write $\lambda_1=n-j$. As in the proof of the upper bound, 
\[
|s_\lambda|\leq\exp\left(-(2+p)\left(\frac jn-\frac{j^2}{n^2}+\frac{j}{n^2}\right)\right).
\]
Thus,
\begin{align*}
\sum_{\lambda:\text{Zone 4-(1)}}d_\lambda|s_\lambda|^k
&=\sum_{j=1}^{\left\lfloor\frac{n}{3}\right\rfloor}\sum_{\lambda_1=n-j}d_\lambda|s_\lambda|^k\\
&\leq\sum_{j=1}^{\left\lfloor\frac{n}{1000}\right\rfloor}\sum_{\lambda_1=n-j}d_\lambda|s_\lambda|^k+\sum_{j=\left\lfloor\frac{n}{1000}\right\rfloor+1}^{\left\lfloor\frac{n}{3}\right\rfloor}\sum_{\lambda_1=n-j}d_\lambda|s_\lambda|^k\\
&=A_1+A_2.
\end{align*}
By (\ref{d*}), 
\begin{align*}
A_1
&\leq\sum_{j=1}^{\left\lfloor\frac{n}{1000}\right\rfloor}\frac{n^j2^{\frac j2}}{\sqrt{j!}}\exp\left(-(2+p)k\left(\frac jn-\frac{j^2}{n^2}+\frac{j}{n^2}\right)\right)\\
&\leq\sum_{j=1}^{\left\lfloor\frac{n}{1000}\right\rfloor}\frac{2^{\frac j2}}{\sqrt{j!}}\exp\left(-(2+p)jc+\frac{j(j-1)(\log n+(2+p)c)+(2+p)j}{n}\right).
\end{align*}
Let $a_j(n)$ denote the $j$-th summand on the right-hand side. Set
\[
b_j(n)=\frac{a_{j+1}(n)}{a_j(n)}=\frac{\sqrt{2}e^{-(2+p)c}}{\sqrt{j+1}}\exp\left(\frac{2j(\log n+(2+p)c)+(2+p)}{n}\right). 
\]
The sequence $b_j(n)$ first decreases and then increases, with a turning point near $j=\dfrac{n}{4(\log n+(2+p)c)}$. Moreover, evaluating the displayed expression at $j=\left\lfloor n/1000\right\rfloor$ gives
\begin{align*}
b_{\left\lfloor\frac{n}{1000}\right\rfloor}(n)
&\leq\frac{\sqrt{2}e^{-(2+p)c}}{\sqrt{\frac{n}{1000}}}\exp\left(\frac{\log n+(2+p)c}{500}+\frac{2+p}{n}\right)\\
&=\sqrt{2000}\exp\left(-(2+p)c+\frac{(2+p)c}{500}+\frac{2+p}{n}\right)n^{\frac{1}{500}-\frac12}\\
&\xrightarrow[n\to\infty]{}0.
\end{align*}
Choose $M$ such that
\begin{equation}\label{lem6.6-2}
\frac{\sqrt{2}e^{-(2+p)c}}{\sqrt{M+1}}\leq\frac14, 
\end{equation}
and take $n$ sufficiently large so that
\[
\exp\left(\frac{2M(\log n+(2+p)c)+(2+p)}{n}\right)\leq2. 
\]
Then
\[
A_1
\leq\sum_{j=1}^{M}a_j(n)+a_M(n)\sum_{i=1}^{\infty}\frac1{2^i}
\xrightarrow[n\to\infty]{}\sum_{j=1}^{M}\frac{2^{\frac j2}}{\sqrt{j!}}e^{-(2+p)jc}+\frac{2^{\frac M2}}{\sqrt{M!}}e^{-(2+p)Mc}.
\]
On the other hand, for $n/1000<j\leq n/3$, we have
\begin{align*}
\sum_{\lambda_1=n-j}d_\lambda|s_\lambda|^k
&\leq\frac{n^j(2e)^{\frac j2}}{j^{\frac j2}}\exp\left(-(2+p)\left(\frac jn-\frac{j^2}{n^2}+\frac j{n^2}\right)k\right)\\
&\leq\frac{(2e)^{\frac j2}}{j^{\frac j2}}\exp\left(\frac{j(j-1)}{n}\log n-(2+p)(cn-1)\left(\frac jn-\frac{j^2}{n^2}+\frac j{n^2}\right)\right)\\
&\leq\exp\left(-\frac j6\log n+\frac j2\log(2000e)+(2+p)|c|(j+1)+\frac{2+p}{3}\right)\\
&\leq\exp\left(-\frac{n}{6000}\log n+O(n)\right). 
\end{align*}
Hence,
\[
A_2\leq\frac n3\exp\left(-\frac{n}{6000}\log n+O(n)\right)\xrightarrow[n\to\infty]{}0.
\]

\noindent\textbf{Zone 4-(2)}

Let $\delta_n=2(1-p)\dfrac{(n+1)(n-3)}{3n^2}$. The argument used for Zone 4-(2) in the proof of the upper bound gives 
\[
\sum_{\lambda:\text{Zone 4-(2)}}d_\lambda|s_\lambda|^k
\leq\left(\sum_{\lambda:\text{Zone 4-(1)}}d_\lambda|s_\lambda|^k\right)(1-\delta_n)^k+\sum_{\substack{\lambda:\text{Zone 4-(2)}\\s_\lambda<0}}d_\lambda|s_\lambda|^k.
\]
For the first term, 
\begin{align*}
&\left(\sum_{\lambda:\text{Zone 4-(1)}}d_\lambda|s_\lambda|^k\right)(1-\delta_n)^k\\
&\leq\left(\sum_{\lambda:\text{Zone 4-(1)}}d_\lambda|s_\lambda|^k\right)\exp\left(-\frac{2(1-p)}{3(2+p)}n\log n+O(n)\right)\\
&\xrightarrow[n\to\infty]{}0.
\end{align*}
For the second term, 
\begin{align*}
\sum_{\substack{\lambda:\text{Zone 4-(2)}\\s_\lambda<0}}d_\lambda|s_\lambda|^k
&\leq\sum_{j=1}^{\left\lfloor\frac{n}{3}\right\rfloor}\frac{n^j2^{\frac j2}}{\sqrt{j!}}\left((1-p)\left(1-\frac{2j}{n}+\frac{2j^2}{n^2}\right)\right)^k\\
&\leq\frac n3\exp\left(\frac{1}{2+p}\log(1-p)\cdot n\log n+O(n)\right)\\
&\xrightarrow[n\to\infty]{}0.
\end{align*}
Consequently,
\[
\sum_{\lambda:\text{Zone 4-(2)}}d_\lambda|s_\lambda|^k\xrightarrow[n\to\infty]{}0.
\]

\noindent\textbf{Zone 4-(3)}

As in the proof of the upper bound, there exists a constant $q_p\in(0,1)$ such that $|s_{(1^n)}|\leq q_p$ for all sufficiently large $n$. Therefore, 
\[
d_{(1^n)}|s_{(1^n)}|^k\leq q_p^k\xrightarrow[n\to\infty]{}0.
\]
\vskip \baselineskip

Combining the estimates above, for every $\epsilon>0$, we choose $M=M(c,\epsilon)\geq1$ such that both (\ref{lem6.6-2}) and
\[
4\cdot\frac{2^{\frac M2}}{\sqrt{M!}}e^{-(2+p)Mc}<\frac{\epsilon}{2}
\]
hold. Then, for all sufficiently large $n$,
\begin{align*}
&\sum_{\lambda_1\leq n-M}d_\lambda|s_\lambda|^k\\
&\leq\left(\sum_{\lambda:\text{Zone 1}}+\sum_{\lambda:\text{Zone 2}}+\sum_{\lambda:\text{Zone 3}}+\sum_{\lambda:\text{Zone 4-(2)}}+\sum_{\lambda:\text{Zone4-(3)}}\right)d_\lambda|s_\lambda|^k+\sum_{j=M}^{\left\lfloor\frac{n}{1000}\right\rfloor}a_j(n)+A_2\\
&\leq\sum_{j=M}^{\left\lfloor\frac{n}{1000}\right\rfloor}a_j(n)+\frac{\epsilon}{2}\\
&\leq4\cdot\frac{2^{\frac M2}}{\sqrt{M!}}e^{-(2+p)Mc}+\frac{\epsilon}{2}\\
&<\epsilon. 
\end{align*}
The desired result follows from this.
\end{proof}

By adapting the arguments used for Lemmas~\ref{lem5.3} and~\ref{lem5.4}, Proposition~\ref{prop5.5}, and Lemma~\ref{lem5.6}, we obtain the following results for the biased random transposition-3-cycle shuffle.

\begin{lem}[{cf.~\cite[Lemma 4.2]{paper2020}}]\label{lem6.3}
Let $l$ be a fixed positive integer. Then, as $n\to\infty$,
\[
\frac{1}{n!}\sum_{\sigma\in S_n}\left|\sum_{j=1}^l\sum_{\lambda_1=n-j}d_\lambda s_\lambda^k\ch^\lambda(\sigma)\right|=\frac{1}{n!}\sum_{\sigma\in S_n}\left|\sum_{j=1}^le^{-(2+p)jc}\Phi_j(\operatorname{Fix}(\sigma))\right|+o(1). 
\]
\end{lem}

\begin{lem}[{cf.~\cite[Lemma 4.4]{paper2020}}]\label{lem6.4}
Let $\epsilon>0$. There exists $M_0=M_0(\epsilon,c)$ such that, for every $M\geq M_0$ and every positive integer $n$,
\[
\left|\frac{1}{n!}\sum_{\sigma\in S_n}\left|\sum_{j=1}^Me^{-(2+p)jc}\Phi_j(\operatorname{Fix}(\sigma))\right|-\frac{1}{n!}\sum_{\sigma\in S_n}\left|\sum_{j=1}^\infty e^{-(2+p)jc}\Phi_j(\operatorname{Fix}(\sigma))\right|\right|\leq\epsilon.
\]
\end{lem}

Define the function $g_c$ by
\[
g_c(x)=e^{-e^{-(2+p)c}}(1+e^{-(2+p)c})^x-1, \qquad x\in\R.
\]

\begin{prop}[{cf.~\cite[Proposition 4.5]{paper2020}}]\label{prop6.5}
For a positive integer $N$,
\[
\sum_{j=1}^\infty e^{-(2+p)jc}\Phi_j(N)=g_c(N).
\]
\end{prop}

\begin{lem}[{cf.~\cite[Lemma 4.6]{paper2020}}]\label{lem6.6}
As $n\to\infty$,
\[
\frac{1}{n!}\sum_{\sigma\in S_n}|g_c(\operatorname{Fix}(\sigma))|\to\E[|g_c(\operatorname{Poiss}(1))|]. 
\]
\end{lem}

\begin{proof}[\normalfont\bfseries Proof of Theorem \ref{thm:biased limit-profile}]
The proof proceeds similarly to that of Theorem~\ref{thm:limit-profile}, by using Lemmas~\ref{lem6.1}, \ref{lem6.2}, \ref{lem6.3}, \ref{lem6.4}, \ref{lem6.6} and Proposition~\ref{prop6.5} in place of Lemmas~\ref{lem5.1}, \ref{lem5.2}, \ref{lem5.3}, \ref{lem5.4}, \ref{lem5.6} and Proposition~\ref{prop5.5}. 
\end{proof}
\vskip \baselineskip

\end{document}